\documentclass{amsart}
\usepackage{amsmath,amsthm,amscd,amssymb,amsfonts}
\ifx\pdfoutput\undefined \usepackage[hypertex]{hyperref} \else \usepackage[pdftex,pdfstartview=FitH]{hyperref} \fi

\DeclareFontFamily{OT1}{rsfs}{}
\DeclareFontShape{OT1}{rsfs}{n}{it}{<-> rsfs10}{}
\DeclareMathAlphabet{\mathscr}{OT1}{rsfs}{n}{it}
\numberwithin{equation}{section}
\newcommand{\absolute}[1]{{\left|{#1}\right|}}
\newcommand{\norm}[1]{{\left\|{#1}\right\|}}
\newcommand{\Satake}{\mathcal{S}}
\newcommand{\Qb}{\mathbf{Q}}

\newcommand{\Lie}{\mathrm{Lie}}
\newcommand{\Orth}{\mathbb{O}}

\newcommand{\temp}{\mathrm{temp}}
\newcommand{\planch}{\mu_{\mathrm{Planch}}}

\newcommand{\haar}{\mu_{\mathrm{Haar}}}
\newcommand{\eigs}{\mathcal{E}}
\newcommand{\R}{\mathbb{R}}
\newcommand{\vol}{\mathrm{vol}}
\newcommand{\Q}{\mathbb{Q}}
\newcommand{\Z}{\mathbb{Z}}

\newcommand{\M}{\mathbf{M}}
\newcommand{\Eis}{E_{1/2+ir}}
\newcommand{\A}{\mathbf{T}}
\newcommand{\B}{\mathbf{B}}
\newcommand{\diam}{{\aleph}}

\newcommand{\Pb}{\mathbf{P}}
\newcommand{\N}{\mathbf{N}}
\newcommand {\equ}[1]     {\eqref{#1}}
\newcommand{\twobytwo}[4] {\begin{pmatrix}
{#1}&{#2}\\{#3}&{#4}\end{pmatrix}}
\DeclareMathOperator{\SL}{SL}

\DeclareMathOperator{\PGL}{PGL}

\newcommand {\C} {{\mathbb C}}
\renewcommand {\H} {{\mathbb H}}
\newcommand{\specnorm}{\mathrm{spec},\infty}

\newcommand{\G}{\mathbf{G}}

\newcommand{\Siegel}{\mathfrak{S}}

\newcommand{\cont}{\mathrm{cts}}
\newcommand{\cusp}{\mathrm{cusp}}
\newtheorem{prop}{Proposition}
\newtheorem{lemma}{Lemma}
\newtheorem{thm}{Theorem}
\newtheorem{cor}{Corollary}
\theoremstyle{remark}

\newcommand{\IGNORE}[1] {}

\usepackage{amsmath,amssymb,amsfonts}
\title {Existence and
Weyl's law for spherical cusp forms}
\author{Elon Lindenstrauss and Akshay Venkatesh}
\begin{document}

\begin{abstract} Let $\G$ be a split adjoint semisimple group over $\Q$ and $K_{\infty}
\subset \G(\mathbb{R})$ a maximal compact subgroup. We shall give
a uniform, short and essentially elementary proof of
the Weyl law for
cusp forms on congruence quotients of $\G(\mathbb{R}) / K_{\infty}$.
This proves a conjecture of Sarnak for $\Q$-split groups,  previously known only for the case $\G=\PGL(n)$.
The key idea amounts to a new type of simple trace formula. 
\end{abstract}

\maketitle

\tableofcontents
\section{Introduction}
Let $M$ be a Riemannian manifold. It was proved by Weyl that if $M$ is
compact, the number of Laplacian eigenvalues less than $T$ is asymptotic to
$c(M) T ^{\dim(M)/2}$. Here $c(M)$ is
the product of the volume of $M$, the volume of the (Euclidean) unit ball in
$\mathbb{R}^{\dim(M)}$, and $(2 \pi)^{-\dim(M)}$.
In general, if $M$ is noncompact but of
finite volume, there is no reason to expect that the Laplacian has any
nontrivial discrete spectrum. In particular, the work of Phillips-Sarnak
\cite{ps} and Wolpert \cite{wolpert} indicates that for a generic nonuniform lattice $\Gamma \subset
\mathrm{PGL}(2,\mathbb{R})$, the Laplacian of the quotient $\Gamma \backslash \mathbb{H}^ 2$ should have only finite discrete spectrum.

On the other hand, Selberg \cite{selberg} has shown that
for $M = \Gamma \backslash \mathbb{H}^ 2$, where $\Gamma$
is a congruence subgroup of $\mathrm{PGL}(2,\mathbb{Z})$,
the Weyl asymptotic holds for the discrete spectrum of the Laplacian,
i.e. the asymptotics of the eigenvalues of the Laplacian behave
as if the surface $M$ were compact. In order to prove this, Selberg
developed his celebrated trace formula.

It has been conjectured by Sarnak 
\cite{Sarnak-cusp-II}
that the same should hold in the general
setting of congruence quotients of noncompact symmetric spaces, that is to
say, $M$ of the form $\Gamma \backslash \G(\R) / K$ where $G$ is a semisimple
algebraic group defined over $\Q$, $K \subset \G(\R)$ a maximal compact
subgroup, and $\Gamma \subset \G(\Q)$ a congruence subgroup. In his thesis, S. D. Miller \cite{miller}
established this conjecture for $\G = \PGL(3)$. More recently, W. M{\"u}ller
\cite{muller} has established this conjecture for $\G = \PGL(n)$. 
The upper bound, i.e. that the number of cusp forms with Laplacian eigenvalue at most $T$ is asymptotically no more than that given by Weyl's bound, was proved for general $\G$ by Donnelly \cite{donnelly}.

The work of Selberg, Miller, and M{\"u}ller
rely essentially on the theory of Eisenstein series, in particular,
estimates on their behavior near the unitary axis. These estimates
-- which are quite delicate in the higher rank case --
are used to explicitly
control the contribution of the noncuspidal spectrum to an
appropriate trace formula (in the work of Miller, to a pretrace formula). 

In the present paper we give a simple proof of the Weyl law for cusp forms,
valid for any split
adjoint group $\G$ over $\Q$, that is to say:
\begin{thm}\label{thm:mainthm}
Let $\G$ be a split adjoint semisimple group over $\Q$, $G_{\infty} = \G(\R)$
and $K_{\infty} \subset G_{\infty}$ a maximal compact subgroup,
$\Gamma \subset \G(\Q)$ a congruence subgroup.
Let $M = \Gamma \backslash G_{\infty}/K_{\infty}$ be the associated
locally symmetric space, endowed with the Riemannian metric that
corresponds to the Killing form on the Lie algebra of $G_{\infty}$
(see page \pageref{BLAH}).  Let $N(T)$ be the number of cuspidal (see
(\ref{eq:cuspdef})) eigenfuntions of the Laplacian with eigenvalue $\leq T$. Then $$N(T) \sim c(M) T ^{\dim(M)/2},$$ where $c(M)$ is as above.
\end{thm}


Our approach is based on the observation that there are strong relations
between the spectrum of the Eisenstein series at different places. This allows us, after passing to an $S$-arithmetic setting, to construct convolution operators whose image can be seen directly (even without knowing anything about Eisenstein series) to be purely cuspidal. In this way, the delicate estimates on Eisenstein series which were the hard part in the works of
Selberg, Miller and M{\"u}ller are entirely avoided. 
Once these operators are constructed (and using Donnelly's upper bound), 
the Weyl law can be proved for general split groups by applying a pretrace formula \equ{eq:pointwise} in a way that is very similar to that used by Miller in his proof of Weyl's law for $\PGL (3, \Z)$.

Already in the classical case of Maass forms on
$\PGL(2,\mathbb{Z})\backslash \PGL(2,\mathbb{R})$ -- i.e., even Maass forms on the modular surface -- the proof seems
to be new. Indeed, we are not aware of any proof that there are any
even cusp forms on $\PGL(2,\mathbb{Z}) \backslash \mathbb{H}$ that
does not require explicit bounds on the constant term of the
Eisenstein series. In order to illustrates several
of the main ideas of the paper in a relatively simple context we give an elementary and self-contained proof for the existence of infinitely many such forms in $\S 2$.

The idea of using convolution operators with purely cuspidal image is not completely new, and is for example the basis of what is known as the ``simple trace formula'": see, for example, \cite{FK}.  However, previous constructions used convolution operators which factor into the composition of independent convolution operators at each place. Such operators have not only the Eisenstein series but any everywhere unramified automorphic form in their kernel (all Hecke-Maass forms on $\PGL(2,\mathbb{Z})\backslash \H$ are by definition everywhere unramified). This type of trace formula was used by Labesse and M\"uller
\cite{labessemuller} to prove a ``weak Weyl law'', but cannot be used
to prove the full Weyl law for precisely the reason just observed. 
 In a sense, the main observation of the present paper
is that there
is a much larger class of convolution operators with cuspidal image than 
those which have been previously utilized.  

We have not aimed for the greatest generality.
The approach described will 
probably apply much more generally
than split adjoint groups over $\Q$; we have restricted ourselves
to this case 
mainly for simplicity of exposition. 
In order both to illustrate some of the techniques and to normalize the constants involved with minimal amount of computation we also prove Weyl's law for compact quotients --- i.e. for groups which have $\Q$-rank zero. For these compact quotients Weyl's law is well known --- indeed, the original Weyl's law is applicable, and moreover an explicit analysis of the error term using the tools of harmonic analysis on Lie groups was given by Duistermaat, Kolk and Varadarajan~\cite{DKV2}.

We do not believe that the present approach supersedes the approach via
the trace formula. For instance, while we have not attempted to extract explicit estimates on the error term, we expect the trace formula approach, when it can be made to work, should at least in theory yield a substantially better error term than what can be extracted from the proof we give here. 
Also, our approach gives no information about cusp forms whose spectral parameters at the different places happen to satisfy the same relations as those satisfied by Eisenstein series.


{\bf Acknowledgements.} 
We would like to thank Erez Lapid, Steve Miller and Peter Sarnak for their encouragement
of this project.
Thanks to Erez Lapid's suggestions, the first proof of Prop. \ref{prop1}
has been improved. Discussions with Erez also clarified to us the limitations of our techniques. Steve Miller provided detailed and very helpful comments on a draft of this paper. We would like to thank Peter Sarnak and Peter Lax for helpful discussions. This research has been conducted while both E.L. and A.V. were Clay Research Fellows; this generous support from the Clay Mathematics Institute is much appreciated. E.L. was also supported in part by NSF grant DMS-0434403; A.V. was supported
in part by NSF grant DMS-0245606.


\section{Existence of cusp forms for $\PGL(2,\Z) \backslash \mathbb{H}$}
\label{section: hyperbolic plane}
We first briefly give a proof of
the infinitude of cusp forms for the special case of $\PGL(2, \Z)
\backslash \H$. In this case, we can interpret
the method quite concretely in terms of a wave equation on hyperbolic space;
in the higher rank case, we will use convolution operators instead. As we will explain at the end of this section, the operator we construct explicitly using the wave equation can also be interpreted as a convolution operator.

Let $ \Delta =- y ^ 2 (\frac {\partial ^ 2 }{ \partial x ^ 2} +
\frac{\partial ^ 2 }{ \partial y ^ 2})$
be the hyperbolic Laplacian on $\H$.
Then $L ^ 2 ( \PGL (2, \Z) \backslash \H) =
\langle 1 \rangle \oplus L ^ 2 _ {\cont} \oplus L ^ 2 _ {\cusp}$,
where $L ^ 2 _ {\cont}$ is spanned by the continuous spectrum of $\Delta$ on
$\PGL (2, \Z) \backslash \H$,
$L ^ 2 _ {\cusp}$ is the cuspidal part, on which $\Delta$ acts discretely.
Explicitly, $L ^ 2_{\cont}$ is spanned by the Eisenstein series
$\Eis$
(see, e.g. \cite{iwaniec}). 
We shall write $L_0 ^ 2$ for the orthogonal complement of the constant
function in $L ^ 2$.

We wish to show that $L ^ 2_{\cusp} \neq \{ 0 \}$.

For any prime $p$ we also have the Hecke operator $T _ p$,
which acts on functions on $\PGL (2, \Z) \backslash \H$ via the rule:
\begin{equation}
\label{definition of Hecke operator}
T_p f(z)= \frac {1 }{ \sqrt p} \left( f (p z) + \sum_ {k = 0} ^ {p - 1} f(
\frac{z+k}{p})\right).
\end{equation}
Then $T_p$ commutes with $\Delta$.
The Eisenstein series $E_{1/2 + ir}(z)$ are joint eigenfunctions of $\Delta$ and
$T_p$:
\begin{align}
\Delta \Eis & =  (\frac {1 }{ 4} + r ^ 2) \Eis \label{eigenvalue of Laplacian} \\
T _ p \Eis & = (p ^ {i r} + p ^ {- i r}) \Eis  \label{eigenvalue of Hecke operator}
.\end{align}

Let us proceed formally for a moment
to indicate the main idea of the method.
From (\ref{eigenvalue of
Laplacian}) and (\ref{eigenvalue of Hecke operator}),
the operator $\diam := T_p - p ^{\sqrt{1/4 - \Delta}}-
p ^{-\sqrt{1/4-\Delta}}$ annihilates
$\Eis$.
The operator $\diam$ may be given a rigorous
interpretation either in terms of the wave equation
or using convolution operators. In the present section
we shall use the wave equation.
In any case, to show that $L ^ 2_{\cusp} \neq \{ 0 \}$
it suffices to find a {\em single nonconstant function}
not annihilated by $\diam$; this
we do by choosing an appropriate test function supported high
in the cusp.  A key ingredient will be the fact
that solutions to the wave equation propagate at finite speed.

We now detail how $\diam$ may be understood
in terms of the wave equation.
Since we have given
a detailed treatment of the (more general) convolution
approach in the subsequent sections,
we will only sketch this approach, omitting careful justification
of issues concerning existence and uniqueness of solutions to the wave
equation.

Equations \equ{eigenvalue of Laplacian} and \equ{eigenvalue of Hecke operator}
admit a nice interpretation in terms of the automorphic wave equation
\begin{equation} \label{automorphic wave equation}
u_{tt} = - \Delta u + \frac {u }{ 4}
.\end{equation}
A solution $u = u(x+iy, t)$ to (\ref{automorphic wave equation})
may be regarded as describing the amplitude of a wave propagating
in the hyperbolic plane.
The low order term of $u/4$ is natural for the hyperbolic Laplacian; see \cite [pp. 7-11]{ Lax-Phillips-automorphic-book}.

Fix $r \in \mathbb{R}$. Consider the unique solution $ u[r] (x + iy, t)$ of \equ{automorphic wave
equation} with initial conditions $u[r]|_{t=0}= \Eis$ and
$(u[r])_t|_{t=0}=0$, i.e. the only solution of \equ{automorphic wave
equation} invariant under the time reversal with the prescribed value
$\Eis$ at
time zero.
Then,
by \equ{eigenvalue of Laplacian}, $u[r] (x + iy, t) = \frac {1 }{ 2}
\Eis (x + iy) (e ^ {itr} + e ^ {- itr})$.
From \equ{eigenvalue of Hecke operator} it follows that:
\begin{equation} \label{propagation of Eisenstein series}
T _ p \Eis(x+iy) = 2 u [ r] (x+iy, \log p)
.\end{equation}

The important property of the above equation is that it gives $T _ p
\Eis$ in terms of propagation by the wave equation for fixed time $\log p$
which does not depend on $r$.

For every $t \in \R$ we can define a linear endomorphism
$U _ t$ of $L ^ 2 (\PGL
(2, \Z) \backslash \H) \cap C ^ \infty (\PGL (2, \Z) \backslash \H)$ to
itself,  taking a function $f(x+iy)$ to $2 u(x+iy,t)$, where $u$ is the
solution to \equ{automorphic wave equation} with $u|_{t=0}=f$,
$u_t|_{t=0}=0$. One may show that this operator is well-defined
in a standard way; moreover, it is self-adjoint (``time reversal
symmetry.'') Formally speaking, one may write
$U_t = e ^{t \sqrt{1/4-\Delta}} + e ^{- t \sqrt{1/4-\Delta}}$;
in fact, $U_t$ gives a rigorous meaning to the right-hand side.

Equation \equ{propagation of Eisenstein series}, and a simple computation for the constant function,  allow us to conclude the following basic fact:

\begin{prop} \label{prev}
For every $f \in L ^ 2 _ {\cont} \oplus \langle 1 \rangle $,
\begin{equation} \label{fundamental identity}
T_p f = U_{\log p} f
.\end{equation}
\end{prop}

Both $T_p$ and $U_{\log p}$ are self-adjoint.
We therefore deduce from Prop. \ref{prev} that:
\begin{cor}
For every smooth $f \in L ^ 2  (\PGL (2, \Z) \backslash \H)$, \ $[T_p -
U_{\log p}] f \in L ^ 2 _ {\cusp}$.
\end{cor}

The operator $T_p - U_{\log(p)}$ thus gives a rigorous interpretation
to $\diam$.

In order to show that $L ^ 2 _ {\cusp}$ is nonzero,
we only need to find one function in $L  ^ 2 _ 0 (\PGL (2, \Z) \backslash \H)$
that does not satisfy \equ{fundamental identity}. We will show that there
are many, by constructing them ``high in the cusp.'' 

For every $R>0$ let $\Omega_R$
be the ``Siegel domain'' $\twobytwo 1 \Z 0 1 \backslash \left\{ x + iy: y>R
\right\}$. One easily sees that for $R >1$, the natural projection 
$\Omega_R
\to \SL (2, \Z) \backslash \H$ is injective,
so the image  of $\Omega_R$ in $\PGL(2, \Z) \backslash \H$ is the quotient of $\Omega_R$ by the reflection
$x+iy \mapsto -x+iy$.  
Let $L^2_{0, even}$ be the subspace of $L^2(\Omega_R)$ consisting of 
functions with integral $0$ and which are invariant under $-x+iy \mapsto x+iy$. 
Then we may regard the
Hilbert space $ L ^ 2 _ {0,even} (\Omega _ R)$ as being embedded isometrically in $L ^ 2 _ 0 (\PGL (2, \Z) \backslash \H)$. Similarly, we define $C^{\infty}_{0, even}(\Omega_R)$. 

For $0 \neq n \in \Z, R >1$ we set $V_{n,R}$ to be the subspace of $C^{\infty}_{0,even}(\Omega_R)$
consisting of functions of the form $f(x+iy) = h(y) \cos(2 \pi nx)$. 

Suppose $R > e ^{t}$. The solutions to \equ{automorphic wave equation}
propagate at speed at most $1$ on the hyperbolic plane.
It follows that for $f \in V_{n,R}$ the function $U_{t} f$
is supported in a $t$-neighbourhood of $\Omega_R$,
i.e. in $\Omega_{R e ^{-t}}$. Moreover, the domain $\Omega_{R e ^{-t}}$
admits an action of $\mathbb{R}/\mathbb{Z}$,
namely $sf(x+iy) = f(x+s+iy)$ for $s \in \mathbb{R}/\mathbb{Z}$.
This action commutes with the Laplacian and so also with $U_t$.
An element $f \in C^{\infty}_{0,even}(\Omega_R)$ belongs to
$V_{n,R}$ exactly when it transforms under some linear combination
of the characters of $s \mapsto e^{2 \pi i n s  }$ and $s \mapsto e^{-2 \pi i n s}$ under this action. 

It follows that $U_t f \in V_{n,R}$ if $f \in V_{n,R}$. We deduce:
\begin{equation} \label{waveprop}
U_t V_{n,R} \subset V_{n, R e ^{-t}} \qquad R > e ^ t.
\end{equation}
On the other hand, it follows from \equ{definition of Hecke operator} that,
for $R > p$, we have:
\begin{equation} \label{heckeprop}
T_p V_{n,R} \subset \begin{cases}
V _ {pn, R/p}& \text{if $p \nmid n$}\\
V _ {pn, R/p} \oplus V_{n/p, pR}& \text{if $p \mid n$}
\end{cases}
.\end{equation}
Since $V_{n,R}$ and $V_{n',R'}$ are orthogonal for every $n \ne n '$,
we conclude from (\ref{waveprop}) and (\ref{heckeprop})
that $T_p V_{n,R} \perp U_{\log(p)} V_{n,R}$ if $n \neq 0$ and $R > p$.
Consequently, if $f \in V_{n,R}$ and $T_p f \neq 0$, then $[T_p -
U_{\log p}]f \neq 0$.

On the other hand, one easily verifies, by looking at the Fourier expansion, that
 $T_p$ is injective on $V_{n,R}$. 

We have proved:

\begin{prop}
For any $n \ne 0$, $R>p$, the
map $(U_{\log(p)} - T_p)$ is an injective map
of $V_{n,R}$ into $L ^ 2_{\cusp} (\PGL(2,\Z) \backslash \mathbb{H})$. 
\end{prop}

Evidently this implies the infinitude of cusp forms
in a quantifiable fashion. Indeed, one may
deduce by a standard variational argument a weak form of the Weyl law,
i.e. that there is a constant $c > 0$ such that
there are at least $c T$ cusp forms of eigenvalue less than $T$.
Since we prove the full Weyl law later, we omit the easy proof.
We refer to the beginning of Sec. \ref{sec:weyllaw}
for an outline of the idea of the proof of the Weyl law.
For now we just comment that one problem is that, because $\diam$ kills the Eisenstein series,
it will also act by a very small scalar on any cusp form that suitably ``mimics'' an Eisenstein series. 
Thus any approach which uses $\diam$ will lose information about these forms; but it will turn out
that  ``fake Eisenstein series'' are spectrally very sparse.

The operator $\diam$ can be interpreted as a convolution operator, but this requires a slightly different viewpoint. Since this will be used later in the proof of the Weyl law, we explain this in some detail.

We first explain what we mean by a convolution operator on $L ^ 2 (\PGL (2, \Z) \backslash \H)$. We identify $L ^ 2 (\PGL (2, \Z) \backslash \H)$ with the subspace  $L ^ 2 (\PGL (2, \Z) \backslash G / K)$ of $K = \mathrm{PO} (2, \R)$-invariant functions on $\PGL (2, \Z) \backslash G $ with $G= \PGL (2, \R)$. If $k \in C ^ \infty _ c (K \backslash G/K)$, i.e. is a bi-$K$-invariant compactly supported smooth function\footnote{Thus $k$ is a function on $G$, not an element of $K$; we hope this notation does not cause confusion.}, it gives rise to an operator $f \to f \star k$ on $L ^ 2 (\PGL (2, \Z) \backslash \H) = L ^ 2 (\PGL (2, \Z) \backslash G / K)$ given by
\begin{equation*}
f \star k (g) = \int_ G f (g h ^{-1}) k (h) d h.
\end{equation*}
More generally, one can consider convolution with compactly supported distributions instead of functions (in this case the convolution operator is well-defined only on suitably smooth functions $f$, for example on $C ^ \infty (\PGL (2, \Z) \backslash \H)$.)

The operator $\diam$ is a sum of two operators: the Hecke operator $T _ p$, and the operator $U _ {\log p} = p ^{\sqrt{1/4 - \Delta}}+
p ^{-\sqrt{1/4-\Delta}}$. The second of these operators, $U _ {\log p}$, is already a convolution operator, though not with a function but with a distribution.

Explicitly, let $\Xi _ s (z)$ denote the spherical function with parameter $s$, i.e. the unique function on $\H$ which is
\begin{enumerate}
\item spherically symmetric, i.e. depending only on the hyperbolic distance of $z$ from the point $i \in \H$
\item satisfies $\Delta \Xi _ s  = (\frac 14 + s ^ 2) \Xi_s$
\item $ \Xi _ s (i) = 1$.
\end{enumerate}
We can identify it with a $K$-bi-invariant function on $G$. 

For any bi-$K$-invariant compactly supported function (or distribution) $k$ and $s$, the function $ \Xi _ s \star k$ also satisfies (1) and (2) above, and so is equal to $ \hat k (s) \Xi _ s$; the map $k \to \hat k (s)$ is called the spectral transform. Its inverse is given by (see \cite[Theorem 4.8] {Helgason-topics}{}\footnote{The hyperbolic metric used by Helgason and the one we use here differ by a factor of 2, and all other quantities defined by the metric, such as the Laplacian, the spherical transform etc. need to be scaled accordingly.})
\begin{equation*}
k(z) = \frac {1 }{ 4 \pi ^ 2} \int_ \R \hat k (s) \Xi _ s (z) \absolute {c(s)} ^{- 2} d s \qquad c(s) = \pi ^ {- \tfrac{1}{2}} \frac {\Gamma (i s) }{ \Gamma (i s+1/2)}
\end{equation*}
Again, we may identify the right $K$-invariant function $k$ with a function on the upper half plane $\mathbb{H}$.

By definition of $U _ {\log p}$, its spectral transform is  $p ^{i s} + p ^{-i s}$, so formally the corresponding spherical distribution is given by
\begin{equation*}
k_{U_{\log p}} (z) = \frac {1 }{ 4 \pi ^ 2} \int (p ^ {i s} + p ^ {- i s}) \Xi _ s (z) \absolute {c(s)} ^ 2 d s.
\end{equation*}
It follows from the finite propagation speed of the wave equation that $ k _ {U _ {\log p }} (z)$ is supported in the hyperbolic disc of radius $\log p$ around $i$; this can also be deduced from a
distributional variant of Paley-Wiener type theorem \cite[Theorem 4.7]{ Helgason-topics}.

One may avoid the use of distributions in the following fashion:
If $H \in C_c ^ \infty (K \backslash G/K)$, then for any $f \in C ^ \infty (\PGL (2, \Z) \backslash \H)$,
\begin{equation*}
U _ {\log p} (f) \star H = f \star k_1
\end{equation*}
with $k_1$ now a smooth compactly supported function given explicitly by
\begin{equation*}
k_1 (z) = \frac {1 }{ 4 \pi ^ 2} \int \hat H(s) (p ^ {i s} + p ^ {- i s}) \Xi _ s (z) \absolute {c(s)} ^ 2 d s
\end{equation*}
(it is easy to see this defines a smooth function; compact support requires use of the properties of the wave equation and/or a Paley-Wiener type theorem).

The operator $T _ p$ cannot be directly interpreted as a convolution operator on $ \PGL (2, \Z) \backslash \H$. However, one has the following fundamental isomorphism: in the notations of the next section, we consider the set of places $S =  \left\{ \infty, p \right\}$, and take $G_S = \PGL (2, \R) \times \PGL (2, \Q _ p)$, $ \Gamma = \PGL (2, \Z [1 / p])$ and $K_S = K \times \PGL (2, \Z _ p)$, with $\Gamma$ considered as a discrete subgroup of $G _ S$ using the diagonal embedding. Then
\begin{equation*}
\Gamma \backslash G_S/ K_S \cong \PGL (2, \Z) \backslash \H.
\end{equation*}
This isomorphism between spaces gives us an isomorphism $ \iota ^ {*}$ of  $L ^ 2 (\PGL (2, \Z) \backslash \H)$ with the space of right $K_S$ invariant functions in $ L ^ 2 (\Gamma \backslash G _ S)$. A straightforward calculation now verifies that for
$f \in L ^ 2 (\PGL (2, \Z) \backslash \H)$
\begin{equation*}
\iota ^*(T _ pf) = \iota ^* (f) \star k_{T_p}
\end{equation*}
for the $K_S$ invariant distribution $k_{T_p}$ on $G_S$ defined by
\begin{equation*}
k_{T_p} (g_\infty, g _ p) = p ^ {- \tfrac{1}{2}} \delta_1 (g_\infty) 1_{K_p \Bigl(\begin{matrix} p &0\\ 0& 1\end{matrix}\Bigr) K_p}\!\!(g_p) \qquad K_p = \PGL (2, \Z _ p),
\end{equation*}
with $\delta_1$ being Dirac's delta measure. 
Again we see that if further convolved by $H_\infty \in C_c ^ \infty (K \backslash G / K)$, or more precisely by $H (g _ \infty, g _ p) = H_\infty (g _ \infty) 1_{K_p} (g _ p)$ (which is now a bi-$K_S$-invariant, compactly supported function on $G_S$) we have that
\begin{equation*}
\iota (T _ p(f) \star H_\infty) = \iota (f) \star k'_1
\end{equation*}
with $k'_1$ the function
\begin{equation*}
k'_1(g_\infty,g_p) = p ^ {- \tfrac{1}{2}} H_\infty (g_\infty) 1_{K_p \Bigl(\begin{matrix} p &0\\ 0& 1\end{matrix}\Bigr)  K_p} \!\!(g_p)
.\end{equation*}

In this way we see that $\diam$ on $L ^ 2 (\PGL (2, \Z) \backslash \H)$ can be viewed as a convolution operator on the space $ L ^ 2 (\Gamma \backslash G _ S/K _ S)$ by the distribution $k_{\diam} =  k_{U_{\log p}} - k_{T_p}$, and for any $H = H_\infty 1_{K_p}$, the composition $f \mapsto (\diam f) \star H$ under the isomorphism $ \iota ^ {*}$ becomes the convolution operator $f \mapsto f \star (k_1 - k'_1)$, with $k_1, k'_1$ continuous compactly supported bi-$K_S$-invariant functions as above.
Convolution operators by smooth functions with a purely cuspidal image will be very useful for us later, when we prove the Weyl law, since they can be given a simple spectral expansion that converges pointwise \equ{eq:pointwise}.

After introducing the necessary formalism in the next section, in \S\ref{sec:conv} we show how to find such convolution operators with cuspidal image for a general group.

\section{The Satake map, Plancherel measure, and the spherical transform} \label{sec:satake}

In this section we review some facts from the harmonic analysis of $S$-algebraic groups, which we will need to generalize the constructions of the previous section and to prove Weyl's law.

Let $S$ be a set of places
of $\Q$ containing $\infty$.  Set $\Q_S = \prod_{v \in S} \Q_v$
and let $\Z[S ^{-1}]$ be the ring of $S$-integers. For $x = (x_v)_{v \in S}\in
\Q_S ^{*}$
we set $|x|_S = \prod_{v \in S} |x_v|_v$.

Let $\G$ be a split adjoint semisimple group over $\Q$; we shall
in fact choose a model defined over $\Z[1/B]$ for some positive integer $B$. 
We shall use bold-face letters to designate algebraic groups
(always over $\Q$) and usual letters to denote their points over
$\R$ or $\Q_S$. By $Z(\star)$ we mean the center of the group $\star$.

We fix (once and for all) a pair $\A \subset \B$
of a $\Q$-split torus $\A$ contained in a Borel subgroup $\B$.
Let $\N$ be the unipotent radical of $\B$, so $\N = R_u(\B)$
and $\B = \A \cdot \N$.
If $\Pb$ is a parabolic subgroup containing $\B$, then
we have a Levi decomposition $\Pb = \M_P
\N_P$.  Let $\A_P = Z(\M_P)$.


Set $G_S = \G(\Q_S), A_S = \A(\Q_S), A_{P,S}
= \A_P(\Q_S), N_S = \N(\Q_S), B_S = \B(\Q_S)$,
and similarly define $P_S, M_{P,S}, N_{P,S}$.
Let $K_{\infty}$ be a maximal compact subgroup of $G_{\infty} := \G(\R)$
with the property that $G_{\infty} = N_{\infty} A_{\infty}^{\circ} K_{\infty}$
is an Iwasawa decomposition, where $A_{\infty}^{\circ}$ is the identity component
of $A_{\infty}$; equivalently, the Cartan involution of $G_{\infty}$ which fixes $K_{\infty}$ should
act by inversion on $A_{\infty}^{\circ}$. 

We assume now that all finite places of $S$ are prime to the integer $B$. 
Set now $K_S = K_{\infty} \cdot \prod_{v \in S, v \neq \infty} \G(\Z_v)$,
where $\Z_v$ is the maximal compact subring of $\Q_v$.
Then $K_S$ is a maximal compact subgroup of $G(\Q_S)$.
We will moreover assume that $S$ has the property that, for each finite $v \in S$
and for each parabolic $\Pb$ containing $\B$, 
$K_v \cap \M_P(\Q_v)$ is the stabilizer in $\M_P(\Q_v)$ of a special
vertex in the building of $\M_P(\Q_v)$; and moreover this vertex belongs to the apartment associated
to the maximal torus $\A$. 
This condition is satisfied for almost all finite $v$, 
as follows from \cite[3.9.1]{Tits-Corvallis}. Moreover, $K_{\infty} \cap \M_P(\R)$
is a maximal compact subgroup of $\M_P(\R)$, being the fixed points
of a Cartan involution; so $K_S \cap M_{P,S}$ is a  maximal compact subgroup of $M_{P,S}$.

\label{BLAH}
We now describe some normalizations of metric and measure. 
Fix the Riemannian metric on the symmetric space $G_{\infty}/K_{\infty}$
so that it corresponds to the Killing form on the Lie algebra
of $G_{\infty}$.  (That is, one identifies the tangent space
to $G_{\infty}/K_{\infty}$ at the identity coset with
the orthogonal complement of $\Lie(K_{\infty})$ in $\Lie(G_{\infty})$;
now use the Killing form to endow it with an inner product.)
This fixes a normalization of the Laplacian operator, and 
also gives a measure on $G_{\infty}/K_{\infty}$.

The map $N_{P,S} \times M_{P,S} \times K_S \rightarrow G_S$
is surjective (Iwasawa).
We equip each $G(\Q_v)$, for $v$ finite, with the Haar measure
which assigns $\G(\Z_v)$ mass $1$. We equip $K_{\infty}$ with the Haar measure
of mass $1$, and then choose the Haar measure on $G_{\infty}$ so that it is compatible
with the measure on $G_{\infty}/K_{\infty}$ arising from the Riemannian metric. 
We arbitrarily choose
Haar measures on the other mentioned groups.

Let $W_S$ be the Weyl group of $A_S$ in $G_S$,
i.e. $\mathrm{Norm}_{G_S}(A_S)/Z_{G_S}(A_S)$.
Let $\Delta$
be a simple system of roots (relative to the system of positive
roots for $\A$ defined by $\N$); thus each $\alpha \in \Delta$
gives a map $\alpha: \A \rightarrow \mathbb{G}_m$.

Let $\delta: A_S \rightarrow \mathbb{R}^{*}$ be the square root of the modular
character of $A_S$ acting on $N_S$, i.e.
$\delta(a) = |\prod_{\alpha \in \Phi^{+}} \alpha(a)|^{1/2}_S$, where $ \Phi^{+}$
is the set of all positive roots for $\A$.

Let $\Gamma$ be a congruence subgroup of $G(\Z[S ^{-1}])$.
It is known that the number of $\Gamma$ orbits
on proper $\Q$-parabolic subgroups of $\G$ is finite.
These orbits (roughly speaking) index the ``cusps'' of $\Gamma \backslash
G_S$.
Fix a set of representatives $\mathcal{R} = \{\Qb_1, \dots, \Qb_r \}$.
For each $1 \leq i \leq I$ choose $\delta_i \in \G(\Q)$ such that
$\delta_i ^{-1} \Qb_{i} \delta_i$ contains $\B$.
Put $\Pb_i = \delta_i ^{-1} \Qb_i \delta_i$,
with unipotent radical $\N_i$ and Levi subgroup $\M_i \supset \A$.

We then put $N_{S,i} = N_{P_i, S}, A_{S,i}=  A_{P_i,S},
M_{S,i} = M_{P_i,S}$.
Further set $\Gamma_i = \delta_i ^{-1} \Gamma \delta_i,
\Gamma_{A,i} = A_{S,i} \cap \Gamma_i, \Gamma_{N,i} = N_{S,i} \cap
\Gamma_i$.  If $S \neq \{ \infty \}$ then it is easy to verify that
each $\Gamma_{A,i}$ is infinite.

%

Let $C ^{\infty}_c(K_S \backslash G_S /K_S)$ be the space
of compactly supported smooth functions on $K_S \backslash G_S /K_S$.

Set $\mathfrak{a}_S = A_{S}/(A_{S} \cap K_S)$. It
is isomorphic to a sum of copies of $\R$ and $\Z$.
Set $\mathfrak{a}_S ^{*}   \stackrel{\mathrm{def}}{=}
\mathrm{Hom}(\mathfrak{a}_S, \mathbb{C}^{*})$.
It is a complex manifold in a natural fashion.
We denote by $\mathfrak{a}_{S,\temp}^{*} \subset
\mathfrak{a}_S ^{*}$ the subset consisting of {\em unitary} characters.
We define the operation of ``conjugation''
$\nu \mapsto \overline{\nu}$ on $\mathfrak{a}_S ^{*}$ via the rule
$\overline{\nu}(a) = \overline{\nu(a)^{-1}}$.
This is an involution that fixes exactly $\mathfrak{a}_{S,\temp}^{*}$.
The set $\mathfrak{a}_{S, \temp}^{*}$ is
an abelian group isomorphic to a product of copies of $\mathbb{R}$
and $\mathbb{R}/\mathbb{Z}$.
\IGNORE{In order to relate to the results of the previous section, we consider the example of $\G= \PGL (2)$ and $S = \left\{ \infty, p \right\}$. In this case, we have that
$\mathfrak{a}_{S, \temp}^{*} \cong \left\{(s, \theta _ p) :s \in \R \text{ and }  \theta \in \R / \Z \right\}$; under this isomorphism, the spherical transform of a function $k \in C ^ \infty _ c (K _ S \backslash G _ S/K _ S)$ satisfies that for any cusp forms $\varphi$ on $\PGL (2, \Z) \backslash \H$ (considered also as a right $K _ S$-invariant function on $\Gamma \backslash G _ S$) with $\Delta \varphi = - (1/4 + s ^ 2)$ and $T _ p \varphi = 2 \cos \theta _ p$
\begin{equation*}
\varphi \star k = \hat k(s, \theta _ p) \varphi
.\end{equation*}
}
Fix a Haar measure $\haar$ on $\mathfrak{a}_{S,\temp}^{*}$.

In the case $S = \{ \infty \}$ we obtain spaces
$\mathfrak{a}_{\infty}$, $\mathfrak{a}_{\infty}^{*}$,
$\mathfrak{a}_{\infty,\temp}^{*}$.
The exponential map gives an isomorphism of
$\mathfrak{a}_{\infty}$ with the Lie algebra $\Lie(A_{\infty})$;
in particular, $\mathfrak{a}_{\infty}$ and $\mathfrak{a}_{\infty,\temp}^{*}$
have the structure of real vector spaces, and $\mathfrak{a}_{\infty}^{*}
\cong \mathfrak{a}_{\infty,\temp}^{*} \otimes_{\R} \mathbb{C}$
the structure of a complex vector space.
We denote by $\langle, \rangle$ the (complex) symmetric bilinear
form on $\mathfrak{a}_{\infty}^{*}$
deduced from the Killing
form on $\Lie(A_{\infty})$.
The form $\langle \cdot, \cdot \rangle$ is positive
definite on $\mathfrak{a}_{\infty,\temp}^{*}$; for
$\nu \in \mathfrak{a}_{\infty,\temp}^{*}$
we put $\| \nu \| = \sqrt{\langle \nu, \nu \rangle}$; we will also use $\norm {\cdot}$ to denote the norm
on $\mathfrak{a}_{\infty}$ coming from the Killing form.

For $\nu \in \mathfrak{a}_S ^{*}$,
we denote by $\nu_{\infty}$ its
image under the natural map $\mathfrak{a}_{S}^{*}
\rightarrow \mathfrak{a}_{\infty}^{*}$.

For $\nu \in \mathfrak{a}_S ^{*}$,
regarding $\nu \delta$ as a character of $A_S$,
extending trivially on $N_S$ and inducing, we obtain
a principal series representation of $G_S$
which has a unique spherical subconstituent $\pi(\nu)$.
The twist by $\delta$ guarantees that this is ``normalized'' i.e. $\pi(\nu)$ is unitary
if $\nu \in \mathfrak{a}_{S,\temp}^{*}$.
For any $\nu$
the conjugate-linear dual
of $\pi(\nu)$
is isomorphic to $\pi(\overline{\nu})$;
in particular, if $\pi(\nu)$ is unitary,
$\pi(\nu)$ must be isomorphic to $\pi(\overline{\nu})$
and $\nu = w \overline{\nu}$ for some $w \in W_S$.

Let $\nu \in \mathfrak{a}_{S}^{*}$,
and fix a $G_S$-invariant sesquilinear pairing $[ \cdot, \cdot]$
bewteen $\pi(\nu)$ and $\pi(\overline{\nu})$.
Let $v ^{0}_{\nu} \in \pi(\nu)$
and $v ^{0}_{\overline{\nu}} \in \pi(\overline{\nu})$ be $K_S$-invariant vectors
such that $[ v ^{0}_{\nu}, v ^{0}_{\overline{\nu}}]
= 1$.
Set $\Xi_{\nu}(g)$ to be the function on $G_S$ defined by
$[ \pi(g) v ^{0}_{\nu}, v ^{0}_{\overline{\nu}} ]$, i.e. $\Xi_{\nu}$ is the
spherical function with parameter $\nu$.

Note that if $\nu_{\infty} \in \mathfrak{a}_{\infty}^{*}$
and $\pi(\nu_{\infty})$ is unitarizable, then
$\langle \nu_{\infty}, \nu_{\infty} \rangle \in \mathbb{R}$.
This follows from the fact that the Casimir operator must operate
on the space of $\pi(\nu_{\infty})$ by a real scalar. In fact,
the Laplacian operator on $G_{\infty}/K_{\infty}$
acts on the spherical vector of $\pi(\nu_{\infty})$ by
a scalar that differs from $\langle \nu_{\infty}, \nu_{\infty} \rangle$
by a fixed (additive) constant.

If $k \in C ^{\infty}_c(K_S \backslash G_S /K_S)$ then for every $\nu \in \mathfrak{a}_S ^{*}$ we have that $\Xi _ \nu$ is an eigenfunction of the corresponding convolution operator. The {\bf spherical transform} $\hat{k}(\nu)$ is defined as the corresponding eigenvalue, i.e.
\begin{equation*}
k \star \Xi _ \nu = \hat k (\nu) \Xi _ \nu
.\end{equation*}

The spherical transform $k \mapsto \hat{k}$
extends to an $L ^ 2$-isometry between $L ^2 (K _ S \backslash G _ S/K _ S)$ and $L ^2 (\mathfrak{a}_S ^{*}, \planch)$ for an appropriate $W_S$-invariant measure $\planch$ on $\mathfrak{a}_S ^{*}$, the
{\bf Plancherel measure}, which is absolutely continuous with respect to $\haar$ (See
\cite{gga} for the archimedean case and
\cite{macdonald} for the $p$-adic case). In particular, for any $k \in C ^{\infty}_c(K_S
\backslash G_S/K_S)$:
\begin{equation}\label{isometry}
\int_{g \in G_S} |k(g)|^ 2 dg = \int_{\nu
\in \mathfrak{a}_{S,\temp}^{*}} |\hat{k}(\nu)|^ 2
d \planch(\nu).\end{equation}

The inverse of the spherical transform is given explicitly in terms of the Plancherel measure by
\begin{equation}\label{eq:planchinversion}
k (g) = \int_ {\mathfrak{a}_{S, \temp}^{*}} \hat k (\nu) \Xi _ \nu (g) d \planch (\nu);
\end{equation}
in particular, we have that
\begin{equation} \label{equation: at the identity}
k(e) =  \int_{\nu
\in \mathfrak{a}_{S,\temp}^{*}} \hat{k}(\nu)
d \planch(\nu).
\end{equation}

Let $d = \dim(G_{\infty}/K_{\infty})$ and $r = \dim(A_{\infty})$. 
We recall the following estimates for the density of the Plancherel measure which follow from the explicit form of the Plancherel measure \cite{gga, macdonald}:
\begin{enumerate}
\item
There is a constant $c_1$ such that, for any positive function $g$
on $\mathfrak{a}_{S,\temp}^{*}$, we have
\begin{equation} \label{eq:planchandhaar}
\int_{\mathfrak{a}_{S,\temp}^{*}} g(\nu) d \planch(\nu)
\leq c_1 \int_{\mathfrak{a}_{S,\temp}^{*}} (1+\| \nu_{\infty}\|)^{d-r} g(\nu)d \haar(\nu).
\end{equation}
\item
Moreover, there exists
$\alpha(\G) > 0 $ such that
\begin{equation} \label{eq:alphadef}
\int_{\nu \in
\mathfrak{a}_{S,\temp}^{*}: \| \nu_{\infty}\| ^ 2 \leq
T}
d \planch(\nu) \sim \alpha(\G) T ^{d/2},\end{equation} as $T \rightarrow
\infty$. \end{enumerate}
\medskip\noindent
Our normalizations imply that $\alpha(\G)$ actually does not depend on the choice of $S$ (subject
to $S$ containing $\infty$).

If $f \in C ^{\infty}_c(K_S \backslash G_S /K_S)$, we define
the {\bf Abel-Satake transform} (also referred to as the Harish transform)
$$\Satake f(a) = \delta(a)^{-1} \int_{n \in N_S}
f(na) dn$$
for $a \in A_S$. $\Satake f$ is $A_S \cap K_S$-invariant;
as such, it may be regarded as a function on $\mathfrak{a}_S$.
The Satake transform is an isomorphism
of $C ^{\infty}_c(K_S \backslash G_S /K_S)$ with
$C ^{\infty}_c(\mathfrak{a}_S)^{W_S}$ (this follows
by combining \cite{warner} in the real case and \cite[Thm 4.1]{cartier} for the $p$-adic case)
and the following diagram commutes:
\begin{equation}
\label{eqn:cd}
\begin{CD}
\ C ^{\infty}_c(K_S \backslash G_S /K_S)
@>k \mapsto \hat{k}>> \mathcal{H}(\mathfrak{a}_S ^{*})\\
@V \Satake VV  @V \mathrm{id} VV\\
C ^{\infty}_c(\mathfrak{a}_S)^{W_S} @>FT>> \mathcal{H}(\mathfrak{a}_S ^{*})
\end{CD}
\end{equation}
Here $\mathcal{H}(\mathfrak{a}_S^{*})$ is the space
of holomorphic functions on $\mathfrak{a}_S^{*}$, the right-hand arrow is the identity map, and
the bottom arrow the (appropriately normalized) Fourier transform on the abelian group
$\mathfrak{a}_S$, i.e., if $f \in C^{\infty}_c(\mathfrak{a}_S)$, then
$FT(f)$, evaluated at $\nu \in \mathfrak{a}_S^*$, equals $\int_{\mathfrak{a}_S} f(a) \nu(a) da$
for a suitably normalized Haar measure $da$. 

Finally, we will make (crucial!) use of the following Paley-Wiener type theorem for the spherical transform on real groups (which is essentially equivalent to the above description of the spherical transform and the properties of the Abel-Satake transform)
\begin{thm} [Gangolli \cite{gangolli}, Helgason \cite{helgason}]\label{theorem: Paley-Wiener}
Suppose $k \in C ^ \infty _ c (K _ \infty \backslash G _ \infty/K _ \infty)$, and suppose that $k (\exp a)=0$ for every $a \in \mathfrak{a}_\infty$ with $\norm a > R$. Then the spherical transform $\hat k (\nu)$ is an entire holomorphic function of $\nu$, invariant under the Weyl group $W_\infty$. Moreover, given any integer $N$ we can find $C_N>0$ such that
\begin{equation} \label{equation: Gangolli}
\absolute {\hat k (\xi + i \eta)} \leq C _ N (1+ (\norm \xi ^2 + \norm \eta ^2) ^ {1/2}) ^ {-N} e ^ {R \norm \xi}, \qquad \xi, \eta \in \mathfrak{a}_{\infty, \temp}^*
.\end{equation}
Conversely, suppose $F$ is a $W _ \infty$-invariant entire holomorphic function of $\nu \in \mathfrak{a}_{\infty}^*$, and $R > 0$ is such that there exists a constant $C _ N$ (for any integer $N$)
so that \equ{equation: Gangolli} holds. Then there exists a unique function $k \in C ^ \infty _ c (K _ \infty \backslash G _ \infty/K _ \infty)$ such that $\hat k = F$. Moreover, $k (\exp a) = 0$ for every $a \in \mathfrak{a}_\infty$ with $\norm a > R$.
\end{thm}

\section{A convolution operator with spherical, cuspidal image}
\label{sec:conv}

Notations being as in the previous section, we recall that $\varphi: \Gamma
\backslash G_S \rightarrow \C$ is {\em cuspidal} if
\begin{equation} \label{eq:cuspdef} \int_{\Gamma_{N,i} \backslash
N_{S,i}} \varphi(\delta_i ng) dn = 0\end{equation}
for almost all $g \in G_S$
and for each $1 \leq i \leq I$.

\begin{prop} \label{prop1}
Suppose $k \in C ^{\infty}_c(K_S \backslash G_S / K_S)$
satisfies one of the (equivalent) conditions:
\begin{equation} \label{eqn:cond}
\sum_{\tau \in \Gamma_{A,i}} \Satake k(\tau a) = 0\, 
\mbox{ whenever } 1 \leq i \leq I, a \in \mathfrak{a}_S \end{equation}
or
\begin{equation} \label{eq:10bis}
\hat{k}(\nu) = 0 \mbox{ whenever } \nu|\Gamma_{A,i} = 0, \,
1 \leq i \leq I.\end{equation}

Then the convolution
$f \mapsto f \star k$ maps $L ^ 1_{\mathrm{loc}}(\Gamma \backslash G_S)$
into the space of cuspidal functions.
\end{prop}

Note that in (\ref{eqn:cond}) we regard $\Gamma_{A,i}$
as a subgroup of $\mathfrak{a}_S$ via the natural maps
$\Gamma_{A,i} \subset A_{S,i} \subset A_S \rightarrow \mathfrak{a}_S$.
The equivalence of the two conditions of the Proposition is verified via (\ref{eqn:cd}).
Remark moreover that the series in (\ref{eqn:cond}) is easily verified
to be absolutely and uniformly convergent (indeed it is, locally, a finite
sum).

We shall give two proofs. The first is shorter, and uses
some (entirely elementary) facts about pseudo-Eisenstein series.
In the second we give an explicit and direct
proof from the definition of ``cuspidal.''
The main idea is that Eisenstein series associated to $\Qb_i(\Q_S)$
arise from automorphic data on $\bar{\Gamma} \backslash \mathbf{M}_{\Qb_i}(\Q_S)$,
where $\mathbf{M}_{\Qb_i}$ is the Levi subgroup of $\Qb_i$,
and $\bar{\Gamma}$ the projection of $\Gamma \cap \Qb_i(\Q_S)$ to
$\mathbf{M}_{\Qb_i}(\Q_S)$. If $\bar{\Gamma}$ has a large center, this
constrains the automorphic data and 
forces relations between the parameters of the Eisenstein series at different places.

{\em First proof:}
The assumption that $\hat{k}(\nu) = 0$ whenever $\nu$ annihilates
$\Gamma_{A,i}$ is easily seen to imply
that $E \star k = 0$, whenever $E$ is a pseudo-Eisenstein
series attached to $\Qb_i$: that is to say,
$E = \sum_{\Gamma \cap \Qb_i(\Q_S) \backslash \Gamma} f(\gamma g)$
for $f \in C ^{\infty}_c(\N_{\Qb_i}(\Q_S) (\Gamma \cap \Qb_i(\Q_S))
\backslash G_S)$, where $\N_{\Qb_i}$ is the unipotent
radical of $\Qb_i$.\footnote{See \cite{mw} for an adelic treatment.}

Let $\check{k}(g) = \overline{k(g ^{-1})}$.
Then, for any $E$ as above and $f \in L ^ 2(\Gamma \backslash G_S)$,
we have $\langle f, E \star k \rangle = \langle f \star \check{k}, E
\rangle$. On the other hand, the pseudo-Eisenstein series attached to proper parabolics span
the orthogonal complement of the cuspidal spectrum
(see \cite[II.1.2]{mw} for a proof in the adelic context, which
easily implies the asserted statement). Therefore $f \star \check{k}$
is cuspidal for any $f \in L ^ 2(\Gamma \backslash G_S)$.

On the other hand (\ref{eqn:cond}) is true for $k$
if and only if it is true for $\check{k}$. Thus
$f \star k$ is cuspidal for $f \in L ^ 2(\Gamma \backslash G_S)$.
Since $k$ is compactly
supported, one may replace $f \in L ^ 2$ by $f \in L ^ 1_{\mathrm{loc}}$.  \qed

{\em Second proof:}
Fix $1 \leq i \leq I$ and take $f \in L ^ 1_{\mathrm{loc}}(\Gamma \backslash
G_S)$. Put $f_i (g)= f(\delta_i g)$; then $f_i
\in L ^ 1_{\mathrm{loc}}(\Gamma_i \backslash G_S)$;
moreover $f_i \star k = (f \star k)(\delta_i g)$.
We wish to show that, under the condition (\ref{eqn:cond}),
we have
\begin{equation}\label{eq:doofus}
\int_{\Gamma_{N,i} \backslash N_{S,i}} (f_i \star k)(ng) = 0\end{equation}
for all $1 \leq i \leq I$ and almost all $g \in G_S$.

Set $K(x,y) = \sum_{\gamma \in \Gamma_i} k(x ^{-1} \gamma y)$
(for $x,y$ lying in any fixed compact set, the
$\gamma$-sum is finite). Then
$(f_i \star k)(x) = \int_{y \in \Gamma_i \backslash G_S}  K(x,y)  f_i(y) dy$.
It follows that (\ref{eq:doofus}) will be satisfied if
for any $x,y \in G$ we have $\int_{\Gamma_{N,i} \backslash N_{S,i}}
K(nx,y) dn = 0$. Equivalently:
\begin{equation}
\sum_{\gamma \in \Gamma_i} \int_{n \in \Gamma_{N,i} \backslash N_{S,i}}
k(x ^{-1} n ^{-1} \gamma  y) dn   = \sum_{\gamma' \in  \Gamma_{N,i} \backslash
\Gamma_i}
\int_{n \in N_{S,i}} k(x ^{-1} n ^{-1} \gamma'  y) dn = 0 \end{equation}
Writing the $\gamma'$-sum as a union over $\Gamma_{A,i}$-cosets we
see it suffices to check that (for every $1 \leq i \leq I,x,y$)
$$\sum_{\tau \in \Gamma_{A,i}} \int_{n \in N_{S,i}} k(x ^{-1} n \tau y) dn =0$$
Using the Iwasawa decomposition,
the bi-$K_S$-invariance of $k$,
and the fact that any $\tau \in \Gamma_{A,i}$
centralizes $M_{S,i}$
and normalizes (in a measure-preserving fashion) $N_{S,i}$,
we see that it suffices to check for $m \in M_{S,i}$:
\begin{equation} \label{eq:intermediate}
\sum_{\tau \in \Gamma_{A,i}} \int_{n \in N_{S,i}} k (n \tau m) dn
=0\end{equation}

It is not hard to verify that the sum and integral
(\ref{eq:intermediate}) converge absolutely; in particular,
we can freely interchange order of summation and integration.

Consider the function $k_i: m \mapsto \int_{n \in N_{S,i}} k(nm) dn$
on $M_{S,i}$. It is a compactly supported
smooth spherical function on $M_{S,i}$, i.e.
$k_i \in C ^{\infty}_c(K_S \cap M_{S,i}\backslash M_{S,i}
/K_S \cap M_{S,i})$. (\ref{eq:intermediate}) amounts to the assumption that
$k_i' \stackrel{\mathrm{def}}{=} \sum_{\tau \in \Gamma_{A,i}} k_i(\tau m)$
identically vanishes. Note that $k_i'$ is a spherical function
on $M_{S,i}$, but no longer compactly supported.

Let $N_{M,i} = N_{S} \cap M_{S,i}$, i.e. the $\Q_S$-points
of the unipotent radical of the Borel subgroup $\B \cap \M_{i}$ of $\M_i$.
To check that $k_i' \equiv 0$,
we use the injectivity of the Satake map for $M_{S,i}$ rather than $G_{S}$.

We proceed formally: no convergence problems arise
from the fact that $k_i'$ is not compactly supported,
since it fails to be compactly supported ``only in central directions.''
In particular the integrals we write down below are all absolutely and
uniformly convergent.

It suffices to check that for all $a \in
A_S$ that $\int_{n_2 \in N_{M,i}} k_i'(n_2 a) = 0$, i.e.
$$\int_{n_2 \in N_{M,i}}
\sum_{\tau \in \Gamma_{A,i}} \int_{n \in N_{S,i}}  k(n \tau n_2 a) = 0.$$
Since $\tau$ commutes with each $n_2$,  we may write this
$$\sum_{\tau \in \Gamma_{A,i}} \int_{n \in N_{S,i}} \int_{n_2 \in N_{M,i}}
k(n n_2 \tau a) dn dn_2 = 0.$$ Since $N_{S} = N_{S,i} \cdot N_{M,i}$
(with a corresponding measure decomposition), this is equivalent to:
$\sum_{\tau \in \Gamma_{A,i}} \Satake k (\tau a) = 0 $. \qed

\begin{lemma} \label{lem2}
Let $A$ be an abelian group isomorphic to a sum
of copies of $\mathbb{R}$ and $\mathbb{Z}$,
let $W$ be a finite group of automorphisms of $A$,
and let
$T = \{ A_1, \dots, A_n \}$ be a
finite collection of nontrivial subgroups
(i.e. $A_j \neq \{ 0 \}$).
Then there exists a nonzero distribution on $A$
which is a finite linear combination of point masses, is $W$-invariant,
and vanishes when evaluated on any $A_i$-invariant function. 

\end{lemma}
\proof
Enlarge $T$ to be $W$-invariant. For $h$ any distribution
on $A$ set $\check{h}(x) = \overline{h(-x)}$.
For $1 \leq j \leq n$, choose an element
$0 \neq a_j \in A_j$ and
let $f_j$ be the compactly supported distribution on $A$
given by $f_j(g) = g(a_j) - g(0)$, for $g \in C ^{\infty}_c(A)$.
Set $f_A = f_1 \star f_2 \star \dots \star f_n$,
$f_B = f_A \star \check{f_A}$, $f =
\sum_{w \in W} f_B ^ w$. It is clear that $f$ is a $W$-invariant
compactly supported distribution on $A$ whose $A_j$-averages
vanish (more correctly, $f$ vanishes when tested on any $A_j$-invariant
function).
To see
that $f \neq 0$, remark that the ``Fourier transform''
of $f_B$ is non-negative and not identically vanishing.
Finally, it is clear by construction that $f$ is a finite linear combination of point masses. 
  \qed

We remark that this Lemma shows in fact that there is a smooth nonzero
$W$-invariant function $f \in C ^{\infty}_c(A)^ W$
such that each of the averages $\int_{x \in A_i}
f(y+x) dx= 0$ (for any $y \in A$). To see this, simply convolve with a smooth $W$-invariant
function on $A$ of small support. 

\begin{cor}\label{cor1}
 There exist nonzero $k \in C^{\infty}_c(K_S \backslash G_S/K_S)$
that satisfy (\ref{eqn:cond}). 
\end{cor}
\proof
Apply the previous remark with $A = \mathfrak{a}_S$
and $T = \{\Gamma_{A,i}\}_{1 \leq i \leq I}$, and use the surjectivity of the Satake transform. \qed

\section{The pre-trace formula and a proof of the Weyl law in the compact case} \label{sec: Weyl law in the compact case}


We put $r = \mathrm{rank}(\G) = \mathrm{dim}(A_{\infty})$ and
$d = \dim(G_{\infty}/K_{\infty})$.
Suppose that the spherical transform of $k \in C _ c ^ \infty (K_S \backslash G_S /K_S)$ is nonnegative and that $k$ is such that the corresponding convolution operator has a purely cuspidal image.

There is a basis of cusp forms on $L^2(\Gamma \backslash G_S/K_S)$ that consists
of eigenfunctions for the convolution algebra $C^{\infty}_c(K_S \backslash G_S/K_S)$. 
 Let $\eigs \subset \mathfrak{a}_S ^{*}$ be the multiset
of eigenvalues of cusp forms on $\Gamma \backslash G_S/K_S$,
counted with multiplicity.
For each $\nu \in \eigs$,
let $\varphi_{\nu}$ be the corresponding eigenfunction.

Put $K(x,y) = \sum_{\gamma \in \Gamma} k(x ^{-1} \gamma y)$. By our condition on cuspidality of $k$ we get the spectral expansion
\begin{equation} \label{eq:spectral}
K(x,y) = \sum_{\nu \in \eigs} \hat{k}(\nu)
\varphi_{\nu}(x)
\overline{\varphi_{\nu}(y)}, \ \ \ (x,y \in G_S). \end{equation}
(\ref{eq:spectral}) is, {\em a priori},
an equality in $L ^ 2$, but it is easy to see that
one can specialize pointwise to obtain
\begin{equation} \label{eq:pointwise}
\sum_{\gamma \in \Gamma} k (g ^{-1} \gamma g) = \sum_{\nu \in \eigs} \hat{k}(\nu)
|\varphi_{\nu}(g)|^ 2,  \ \ \ (g \in G_S) .\end{equation}

Let $\Omega \subset \Gamma \backslash G_S$ be any compact set (which we will take to have almost full measure).
Let $\tilde{\Omega}$ be a compact subset of $G_S$
whose projection contains $\Omega$ so that the projection is almost everywhere one to one. Integrating \equ{eq:pointwise} over $\tilde \Omega$ we have
\begin{equation} \label{pretrace formula}
\begin{aligned}
k (1) \vol (\Omega) +
\sum_{\gamma \in Z} \int_ {\tilde \Omega}  k (g ^{-1} \gamma g) dg & = \int_ {\tilde \Omega} \sum_{\gamma \in \Gamma} k (g ^{-1} \gamma g) dg   \\
& = \sum_{\nu \in \eigs} \hat{k}(\nu) \int_ {\tilde \Omega}
|\varphi_{\nu}(g)|^ 2 dg \leq \sum_{\nu \in \eigs} \hat{k}(\nu)
\end{aligned}
\end{equation}
where
\begin{equation} \label{equation defining Z}
Z = (\Gamma \setminus \left\{ e \right\}) \cap \left\{ g x g ^{-1} : g \in \tilde \Omega, \text{and $x$ in the support of $k$} \right\};
\end{equation}
in particular, $Z$ is finite and depends only on $\tilde{\Omega}$ and the support  of $k$.

The idea of using integrating over a compact $\Omega$ and obtaining an inequality as in \equ{pretrace formula}, rather than integrating over all of $\Gamma \backslash G$ and obtaining the usual trace formula, 
was used in Miller's proof \cite{miller} of the Weyl law for cusp forms in $\PGL(3)$.  Our
proof will follow the same pattern, but we will be able to avoid entirely the Eisenstein series. 

\subsection {Weyl's law for compact quotients} \label{sec:Weyllawcompact}
We now briefly describe how to prove Weyl's law for compact quotients, starting
with \equ{pretrace formula}. Of course, this was
proven by Weyl, and it is also a trivial consequence of the usual trace formula.
Nevertheless we take the opportunity to explain
it in the notation of the present paper, as a
precusor to the more involved case of noncompact quotient.

So far in our discussion in this section we have not assumed anything on $\G$, and in particular everything we have said remains valid for compact quotient $ \Gamma \backslash G _ S$. Of course in this case $\G$ is not $\Q$-split, but this was only used for constructing convolution operators with cuspidal image. If $\Gamma \backslash G _ S$ is compact all of the spectrum is cuspidal, and so any convolution kernel in $C_c ^ \infty (K _ S \backslash G _ S/K _ S)$ with positive spherical transform will satisfy \equ{pretrace formula} (indeed, in this case we may as well take $\Omega = \Gamma \backslash G _ S$ in which case \equ{pretrace formula} becomes an identity).

For this subsection, we shall vary the notation slightly and assume that $\G$ is an $\Q$-anisotropic semisimple $\Q$-group, 
$S$ an arbitrary set of places containing $\infty$, $G = \G(\Q_S)$, and $\Gamma$ is a congruence subgroup of $\G(\Z[S ^{-1}])$, so that $\Gamma \backslash G _ S$ is compact.
Since we are interested in this case mostly as an introduction to the split case, we will further assume that $\Gamma$ is torsion free, which avoids some minor difficulties.

To obtain Weyl's law for compact quotients, we apply \equ{pretrace formula}, with
$\Omega = \Gamma \backslash G_S$,  for a family of test functions which will be defined using a function $H ^ \infty \in C ^ \infty _ c (K _ \infty \backslash G _ \infty/K _ \infty)$. This function will be chosen to have the following properties:

\begin{lemma} \label{lem:hchoice}
Let $0<\varepsilon <1$.
There exists $H ^{\infty} \in C ^{\infty}_c(K_{\infty} \backslash
G_{\infty}/K_{\infty})$ so that its spherical transform $h(\nu) := \widehat{H ^{\infty}}(\nu)$ satisfies:
\begin{enumerate}
\item \label{third} $h(\nu)$
depends only on $\langle \nu, \nu \rangle$, for $\nu \in \mathfrak{a}_{\infty}^{*}$.
\item \label{zeroth} $h(\nu)$ is real and non-negative whenever
$\nu$ and $\overline{\nu}$ are $W_{\infty}$-conjugate to each other
(in particular, when $\pi(\nu)$ unitary).
\item \label{first}
$h(\nu) \leq 1$ for every $\nu \in \mathfrak a _ \infty ^{*}$ with
$\langle \nu,\nu \rangle \in \R ^ +$.
\item \label{firstb} $h(\nu) \geq 1-\varepsilon$ whenever $\nu \in \mathfrak{a}^{*}_{\infty}$
and $\langle \nu,\nu \rangle \in \mathbb{R}$, $0 \leq \langle \nu, \nu \rangle \leq 1-\varepsilon$.
\item \label{second} As $t \rightarrow 0$, we have
$$t ^ d
\int_{\mathfrak{a}_{S,\temp}^{*}} h(t \nu_{\infty})
d \planch(\nu) = \alpha(\G)+ O( \varepsilon) ,$$
\item 
\label{decay}
\begin{equation*}
\sup_ {\norm {\nu} \geq 1, \nu \in \mathfrak a _ {\infty,\temp} ^{*}} (1 + \norm {\nu})^{d+1} \absolute { h (\nu)} < \epsilon
\end{equation*}
\end{enumerate}
\end{lemma}

We defer the somewhat technical and not very illuminating proof of this lemma to the end of this section, showing first how, using this $H ^ \infty$, Weyl's law can be proved.

For $0 < t \leq 1$, let $H ^{\infty}_t \in C ^{\infty}_c(K_{\infty}\backslash
G_{\infty}/K_{\infty})$ be so that $\widehat{H ^{\infty}_t}(\nu)
= h(t \nu)$. That this is well defined
is a consequence of the Paley-Wiener type theorem (Thm. \ref{theorem: Paley-Wiener}) for symmetric spaces.
Indeed, Theorem \ref{theorem: Paley-Wiener} 
implies that given any open neighborhood of the identity in $K_{\infty}
\backslash G _ \infty/K_{\infty}$, if $t$ is small enough, the support of $H ^{\infty}_t$ is contained in this neighborhood. We construct a function $H_t \in C ^ \infty _ c (K _ S \backslash G _ S/K _ S)$ from $H_t ^ \infty$ as follows: 
for each finite place $v \in S$ and any $t$, let
$H_t ^ v$ be the characteristic function of $K_v$, and define $H_t
\in C ^{\infty}(K_S \backslash G_S /K_S)$ via $H_t := \prod_{v \in
S} H_t ^ v.$
The spectral transform of $H_t$ is the function
$\nu \mapsto h(t \nu_{\infty})$, where $\nu \mapsto \nu_{\infty}$
is the natural projection $\mathfrak{a}_S ^{*} \rightarrow
\mathfrak{a}_{\infty}^{*}$.

Since $\tilde \Omega$ is fixed, and since the support of $H ^ \infty_t$ shrinks to $K_{\infty}$, for $t$ small enough the set $Z$ defined in \equ{equation defining Z} 
consists only of nontrivial $\gamma \in \Gamma$ which belongs to some conjugate of $K_S$.
The assumption that $\Gamma$ is torsion free shows this set to be empty. 
Now take $t $ to be sufficiently small so that $Z$ is empty. Then
\equ{pretrace formula}, applied with $k = H_t$ and $\Omega = \Gamma \backslash G_S$,  implies that
\begin{align*} \sum_{\nu \in \eigs} h(t \nu_{\infty}) = H_t(1) \vol (\Gamma \backslash G_S)  \end{align*}
By \equ{equation: at the identity} and (\ref{second})  of Lemma~\ref{lem:hchoice},
\begin{equation*}
H _ t (1) = \int_{\nu
\in \mathfrak{a}_{S,\temp}^{*}} h(t \nu _ \infty)
d \planch(\nu) = t ^ {- d} (\alpha (\G) + O(\epsilon))
\end{equation*}
so
\begin{equation} \label{pretrace formula2}
\sum_{\nu \in \eigs} h(t \nu_{\infty}) = t ^ {- d} (\alpha (\G) + O(\epsilon)) \vol(\Gamma \backslash G_S) 
.\end{equation}
Since for any $\nu \in \eigs$, the inner product $\langle \nu _ \infty, \nu _ \infty \rangle = \norm {\nu_{\infty}} ^2 $ differs from the Laplacian eigenvalue of the corresponding eigenfunction $\varphi _ \nu$ by a fixed real constant, we have by \equ{first} of Lemma~\ref{lem:hchoice} that
\begin{equation} \label{equation about h}
h(t \nu _ \infty) \leq 1 \qquad \text{for all but finitely many $\nu \in \eigs$};
\end{equation}
and for all $\nu \in \eigs $ we have that $h(t \nu _ \infty)$ is nonnegative, real, and uniformly bounded for  $t \in (0,1]$.

In view of (\ref{firstb}) of Lemma~\ref{lem:hchoice}, equation \equ{pretrace formula2} implies immediately that for all $t$ sufficiently small
\begin{equation} \label{upper bound}
\absolute {\left\{ \nu \in \eigs: \norm \nu _ \infty < t ^{-1}  \right\}} \leq t ^ {- d} (\alpha (\G) + O(\epsilon))
\vol(\Gamma \backslash G_S) 
.\end{equation}
On the other hand, applying \equ{upper bound} together with Lem. \ref{lem:hchoice}, (\ref{decay}) and \equ{equation about h} gives that
\begin{equation} \label{tail estimate}
\begin{aligned}
\absolute {\sum_ {\nu \in \eigs: \norm {\nu _ \infty} > t ^{-1}} h (t \nu _ \infty)}& \leq \sum_ {n=0}^{\infty} \absolute {\sum_ {\nu \in \eigs: \norm {t \nu _ \infty} \in [2^n, 2^{n+1}] } h (t \nu _ \infty)} \\
& \leq \sum_ n \frac {2 \epsilon  t ^ {- d} \alpha (\G)  \vol(\Gamma \backslash G_S) 2^{(n+1)d} }{ (1 + 2^{n})^{d+1} }= O (\epsilon t^{-d})
\end{aligned}
\end{equation}
so by \equ{pretrace formula2} and \equ{equation about h}
\begin{equation} \label{eq:spongebob}
\begin{aligned}
\absolute {\left\{ \nu \in \eigs: \norm \nu _ \infty < t ^{-1}  \right\}} & \geq O(1) + \sum_{\nu \in \eigs: \norm {\nu _ \infty} < t ^{-1}} h(t \nu_{\infty})
\\
& \geq  \sum_{\nu \in \eigs} h(t \nu_{\infty}) - O(\epsilon t^{-d}) =  t ^ {- d} (\alpha (\G) - O (\epsilon)) \vol(\Gamma \backslash G_S),
\end{aligned}
\end{equation}
concluding the proof of Weyl's law in the compact case.

\proof [Proof of Lemma~\ref{lem:hchoice}]
Let $\Orth$ be the orthogonal group of
$(\mathfrak{a}_{\infty,\temp}^{*}, \langle \cdot, \cdot
\rangle)$.
Then $\Orth$ is a maximal compact subgroup of its complexification
$\Orth_{\C}$, which is identified with the orthogonal
group of $(\mathfrak{a}_{\infty}^{*}, \langle \cdot, \cdot \rangle)$.
Moreover, $W_{\infty} \subset \Orth$.

Let $\chi$ be the characteristic function of the ball of radius $1$ in $\mathfrak{a}_{\infty,\temp}^*$. 
There exists a nonempty open set of Schwarz functions $\psi$ on the real Euclidean
space $\mathfrak{a}_{\infty,\temp}^{*}$ satisfying:
\begin{eqnarray}
\label{eq:k}0 \leq \psi(\nu) < 1, \ \ \ (\|\nu\| \leq 1, \nu \in
\mathfrak{a}_{\infty,\temp}^{*}) \\
\label{eq:ki} \psi(\nu) > \sqrt{1-\varepsilon}, \ \ \ (\|\nu\| \leq \sqrt{1-\varepsilon},
\nu \in \mathfrak{a}_{\infty,\temp}^{*}) \\
\label{eq:l} \sup_{\|\nu\| \geq 1, \nu \in \mathfrak{a}_{\infty,\temp}^{*}}
(1+\|\nu\|)^{d+1} |\psi(\nu)| < \varepsilon/2 \\
\label{eq:m}
\int_{\mathfrak{a}_{\infty,\temp}^{*}}
\left|\psi(\nu) - \chi(\nu)\right| (1+\|\nu\|)^{d-r} d\nu <
\varepsilon/2.
\end{eqnarray}
Since the set of functions whose Fourier transform has compact support is dense in the Schwarz space, we can find a $\psi$ satisfying \equ{eq:k}--\equ{eq:m} whose Fourier transform has compact support.
Averaging over $\Orth$ does not affect the validity of \equ{eq:k} -- \equ{eq:m}; we may thereby assume that:
\begin{equation}\label{eq:inv} \psi(k \nu) = \psi(\nu), \ \ \ \nu \in
\mathfrak{a}_{\infty,\temp}^{*}, k \in \Orth.\end{equation}
Since $W_{\infty} \subset \Orth$, $\psi$ is $W_{\infty}$-invariant.
Moreover, since the Fourier transform of $\psi$ is compactly
supported, it follows that $\psi$ extends to a holomorphic
function on $\mathfrak{a}_{\infty}^{*}$.
It follows that, considered as a function on $\mathfrak{a}_{\infty}^{*}$,
$\psi$ is actually $\Orth_{\C}$-invariant.
Since two nonzero vectors $\nu_1, \nu_2 \in \mathfrak{a}_{\infty}^{*}$ are
$\Orth_{\C}$-conjugate precisely when $\langle \nu_1, \nu_1 \rangle
= \langle \nu_2, \nu_2 \rangle$,
we deduce that $\psi(\nu)$ depends only on $\langle \nu,\nu \rangle$
for $\nu \in \mathfrak{a}_{\infty}^{*}$.
Note that, for $\varepsilon$ sufficiently small,
the conditions (\ref{eq:k}) and (\ref{eq:l}) guarantee
that $\sup_{\nu \in \mathfrak{a}_{\infty,\temp}^{*}} |\psi(\nu)| \leq 1$.

By Theorem~\ref{theorem: Paley-Wiener}, there is a function $F \in C ^{\infty}_c(K_{\infty} \backslash G_{\infty} /
K_{\infty})$ whose spherical transform $\hat F = \psi$. Put $H ^{\infty} = F \star \check{F}$,
where $\check{F}(g) := \overline{F(g ^{-1})}$.

Then $h(\nu) = \psi(\nu) \overline{\psi(\overline {\nu})}$,
whence claim (\ref{zeroth}) follows at once.
Since, for $\nu \in \mathfrak{a}_{\infty,\temp}^{*}$,
we have in fact $h(\nu) = |\psi(\nu)|^ 2$,
we obtain also
claim (\ref{third}), (\ref{first}), (\ref{firstb}) and (\ref{decay}).

Let $c_1$ be as in (\ref{eq:planchandhaar}).
Then:
\begin{multline} \label{ghazal}
\limsup_{t \rightarrow 0} \left|t ^{d} \int_{\mathfrak{a}_{S,\temp}^{*}} h(t \nu_{\infty}) d \planch(\nu) - \alpha(\G) \right|
\\ \leq c_1 \limsup_{t \rightarrow 0} \int_{\mathfrak{a}_{S,\temp}^{*}}
|t ^ d h(t \nu_{\infty}) - t ^ d \chi(t \nu_{\infty})|
(1+\| \nu_{\infty}\|)^{d-r} d \haar(\nu)
\\ \leq c_1 \limsup_{t \rightarrow 0} \int_{\mathfrak{a}_{S,\temp}^{*}}|h(\nu_{\infty})-\chi(\nu_{\infty})|
(t + \| \nu_{\infty}\|)^{d-r} d \haar(\nu)
\\ \leq c_1 \int_{\mathfrak{a}_{S,\temp}^{*}}
|h(\nu_{\infty})-\chi(\nu_{\infty})| (1+\| \nu_{\infty}\|)^{d-r}
d \haar(\nu).
\end{multline}

Recalling that $\sup_{\nu \in \mathfrak{a}_{\infty,\temp}^{*}} |\psi(\nu)| \leq 1$
we see that 
$|h(\nu) - \chi(\nu)| \leq 2 |\psi(\nu) - \chi(\nu)|$ for $\nu \in
\mathfrak{a}_{\infty,\temp}^{*}$.
Applying (\ref{ghazal}) and (\ref{eq:m}), we see that
there is $c_2 > 0 $ such that $$\limsup_{t \rightarrow 0}
\left| t ^ d
\int_{\nu \in \mathfrak{a}_{S,\temp}^{*}} h(t \nu_{\infty}) d \planch(\nu) - \alpha(\G) \right| \leq  c_2 \varepsilon.$$
Reducing $\varepsilon$ as necessary, we conclude that
$h$ satisfies also claim (\ref{second}) of the Lemma.
\qed

\section{The proof of Weyl law for $\Q$-split groups} \label{sec:weyllaw}

We now return to the more interesting case of $\G$ a $\Q$-split group, in which case $\Gamma \backslash G _ S$ is not compact, and the existence and the abundance of cusp forms is substantially more delicate because of the presence of continuous spectrum. As we have remarked in the introduction, the other approaches to the Weyl law in this context have required explicit
estimates on the Eisenstein series to show that they are spectrally ``negligible''; we will avoid this 
issue by taking, from the outset, convolution operators that satisfy (\ref{eqn:cond}) and therefore
kill the continuous spectrum. 



Our approach will be similar to the one presented in the previous section. We proceed by expanding spectrally a family of convolution operator with kernels $k_t \in C _ c ^ \infty (K _ S \backslash G _ S / K _ S)$ for $t \to 0$. Actually,
to derive the full Weyl law we will need to use more than one such family. 
We will avoid the complications of having continuous spectrum in the spectral expansion of $f \mapsto k_t \star f$ by taking kernels $k_t$ satisfying \equ{eqn:cond}, and hence the corresponding
convolution operators will have a purely cuspidal image.

It is clear from \equ{eqn:cond} that the support of the kernels $k_t$ cannot possibly shrink to $ K_S$ as in the compact case. However, by using fairly crude estimates we will show that for $t$ small the terms $\gamma \ne e$ in the finite volume analog to
\equ{pretrace formula} are still negligible.

It also follows from the construction of \S \ref{sec:conv} that even if one is interested only in studying cusp forms on $ \Gamma_\infty \backslash G_\infty / K_\infty$ with $\Gamma_\infty < \G(\Q)$ a congruence subgroup, one needs take $S$ to contain at least one finite place. The passage from Weyl law on the $S$-algebraic $ \Gamma \backslash G_S / K_S$ to that on $ \Gamma_\infty \backslash G_\infty / K_\infty$ is fairly straightforward and is explained in \S \ref{subsec:final}. 
In the case of $\PGL (2, \Z) \backslash \H$ this is particularly straightforward: we take for example $S = \left\{ \infty, p \right\}$, $\G= \PGL (2)$ and note that $ \PGL (2, \Z) \backslash \H$ can be identified with $ \G (\Z [1/p]) \backslash G_S / K_S$, and so Weyl's law on $ \G (\Z [1/p]) \backslash G_S / K_S$ is equivalent to that on $\PGL (2, \Z) \backslash \H$.

Our kernels $k_t$ will be of the form $k \star H_t$ with $H_t$ as in the previous section, and $k$ a compactly supported spherical {\em distribution} satisfying an appropriate version of \equ{eqn:cond}.
Again returning to the special case of $\PGL (2, \Z) \backslash \H$, at the end of section \ref{section: hyperbolic plane} we have shown how the operator $\diam$ defined on functions on $\PGL (2, \Z) \backslash \H$ which has purely cuspidal image corresponds to such a spherical distribution on $K_S \backslash G_S / K_S$.

In order to get Weyl's law we will actually need $k$ to be a distribution which approximates in an appropriate way the projection of $L^2 (\Gamma \backslash G_S/K_S)$ on its cuspidal subspace $L ^ 2_{\cusp} (\Gamma \backslash G_S/K_S)$.

For example, for $\PGL (2)$, we can use the following procedure: first note that the spectrum of $\diam$ is contained in $[-4  \sqrt p, 4 \sqrt p]$. By Weierstrass approximation theorem we can find a  polynomial $p$ with real coefficients and no constant term,  so that $p(t) \approx 1$ in the range $ \epsilon < \absolute t \leq 4 \sqrt p$, and take $k$ to be the distribution corresponding to the operator $p(\diam) ^2$.

\subsection{General case: proof of lower bound in Weyl's law}
We'll now prove the lower bound of the Weyl law (the full Weyl law follows
by combining this with the work of Donnelly; see \S \ref{subsec:final} for details). 

Throughout this section we shall assume that, in addition to those conditions
prescribed in Sec. \ref{sec:satake}, $S$ contains at least one finite place
in addition to $\infty$. 
We denote by $\mathcal{D}(K_S \backslash G_S / K_S)$
the space of $K_S$-bi-invariant, compactly supported {\em distributions}
on $G_S$. 
One defines the spherical transform
$\hat{k}(\nu)$ for $k \in \mathcal{D}(K_S \backslash G_S/K_S)$
and $\nu \in \mathfrak{a}_{S}^{*}$ to be the scalar
by which $k$ acts on the spherical vector in $\pi(\nu)$.
For any $k \in \mathcal{D}(K_S \backslash G_S / K_S)$ we put (the ``spectral norm'') 
$\| k \|_{\specnorm} = \sup_{\nu \in \mathfrak{a}_S ^{*}: \pi(\nu) \, \mathrm{unitary}}
|\hat{k}(\nu)|$, if finite.

We remark that our use of distributions is entirely for {\em notational convenience}. Indeed,
the distributions we will use always be convolved with a smooth function when
they are used in the eventual applications. Thus,
although we will at a certain point use an analogue of the Paley-Wiener theorem
Thm. \ref{theorem: Paley-Wiener} for distributions, this is cosmetic and can be replaced
with the use of Thm. \ref{theorem: Paley-Wiener} as stated.

Let $k \in \mathcal{D}(K_S \backslash G_S/K_S)$
satisfy (\ref{eq:10bis}) and be so that $\| k \|_{\specnorm} < \infty$.
Let $H^{\infty} \in C_c^{\infty}(K_S \backslash G_S/K_S)$;
defining, for $t > 0$,  $H_t$ as in Sec. \ref{sec:Weyllawcompact}, 
we put $k_t = k \star H_t$. Then $k_t \in C_c^{\infty}(K_S \backslash G_S/K_S)$
and also satisfies (\ref{eq:10bis}).

The space $\Gamma \backslash G_S/K_S$ may be regarded as a finite union of locally symmetric spaces, each with universal covering $G_{\infty}/K_{\infty}$.
Define $\eigs$ as prior to (\ref{eq:spectral}). 
It follows from the result of Donnelly \cite{donnelly} that
\begin{equation} \label{eq:donnelly}
\#\{ \nu \in \eigs: \langle \nu_{\infty},\nu_{\infty} \rangle \leq T \} \leq C T ^{d/2},\end{equation}
for some constant $C$. (We will have need of the precise constant, but only at the end of the argument:
see 
(\ref{eq:upperbound}) where we specify the constant. For most of the argument, any constant $C$ would do).

\begin{lemma} \label{lem:main}
 Let $F \subset G_S$ be a compact subset
so that the projection of $F$ to $G_{\infty}$
does not intersect $K_{\infty}$.
Let $k \in \mathcal{D}(K_S \backslash G_S/K_S)$ be such that
$\| k \|_{\specnorm} < \infty$.
Then there is a constant $c$, depending possibly on $k$ and $H ^{\infty}$,
such that
\begin{equation}\label{eq:kestclose}
\sup_{g \in F} |k_t(g)| \leq c
t ^{-d+1/2}, \ \ \ (0 < t \leq 1).  \end{equation}

Moreover, there exists a sequence $k_n \in \mathcal{D}(K_S \backslash G_S/K_S)$
such that $\| k_n \|_{\specnorm} < \infty$ and moreover:
\begin{enumerate}
\item \label{mainfirst} Each $k_n$ satisfies (\ref{eq:10bis}),
\item \label{mainsecond} $\widehat{k_n}(\nu) \geq 0$ whenever $\pi(\nu)$ is unitary, and:
\item \label{mainthird}
Suppose $h$ is a non-negative Schwarz function on
$\mathfrak{a}_{\infty,\temp}^{*}$. Then
$$\lim_{n \rightarrow \infty}
\left( \liminf_{t \rightarrow 0} \frac{\int \widehat{k_n}(\nu)
h(t \nu_{\infty})
d \planch(\nu)}{\| k_n \|_{\specnorm}
\int h(t \nu_{\infty}) d \planch(\nu)}\right)= 1.$$
\end{enumerate}
\end{lemma}

We defer the (tedious but straightforward) proof of this Lemma until the next section.

\proof (of the lower bound in the Weyl law)

We shall take $k \in \mathcal{D}(K_S \backslash G_S/K_S)$ to be one of the sequence $k_n$ constructed in Lem. \ref{lem:main},
and $H^{\infty} \in C^{\infty}_c(K_{\infty} \backslash G_{\infty}/K_{\infty})$
will be as provided by Lem. \ref{lem:hchoice}. In particular, 
with these choices,  $\hat{k}(\nu) \geq 0$ and  $h(t \nu_{\infty}) \geq 0$
for all $\nu \in \eigs$ and all $0 \leq t \leq 1$.  Moreover, $k_t$ satisfies (\ref{eq:10bis}). 

Let $\Omega \subset \Gamma \backslash G_S$ be any compact set.
Let $\tilde{\Omega}$ be any compact subset of $G_S$
whose projection to $\Omega$ is almost everywhere one to one.
Define $Z$ as in \equ{equation defining Z}, except replacing
``support of $k$'' by any compact set that contains
the support of $k_t$ by $0 < t < 1$.

For each $\gamma \in Z$, let $F_{\gamma}$ be the set of $g_{\infty} \in G_{\infty}$
such that $g_{\infty}^{-1} \gamma g_{\infty}$ belongs to $K_{\infty}$.
Then $F_{\gamma}$ has zero measure
(note that this requires that $\gamma$ is not central, so uses the adjointness of $\G$).
Let $\tilde{\Omega}' \subset \tilde{\Omega}$ be a compact subset of $\tilde{\Omega}$ whose projection to $G_{\infty}$
does not contain $\cup_{\gamma \in Z} F_{\gamma}$; since each $F_{\gamma}$ is a null set, we may choose $\tilde{\Omega}'$
to have measure arbitrarily close to that of $\tilde{\Omega}$.

From \equ{pretrace formula}, with $k=k_t$,  we deduce:

\begin{equation}\label{eq:dnc}
\sum_{\nu \in \eigs} \hat{k}(\nu)
h(t \nu_{\infty}) \geq
k_t (1) \vol (\tilde{\Omega}') +
\sum_{\gamma \in Z} \int_ {\tilde \Omega'}  k (g ^{-1} \gamma g) dg\end{equation}

The set $\{ g ^{-1} \gamma g: g \in \tilde{\Omega}', \gamma \in Z \}$ is a compact
subset of $G_S$, whose projection to $G_{\infty}$ does not intersect $K_{\infty}$. So, for all $\gamma \in Z$,
we have by  (\ref{eq:kestclose})  that $\left| \int_{g \in \tilde{\Omega}'} k_t(g ^{-1} \gamma g) \right| \leq c t ^{-d+1/2}$,
for $0 \leq t \leq 1$, and for some constant $c$ depending on $k, H ^{\infty}, \tilde{\Omega}$
and $\tilde{\Omega}'$. So, applying Plancherel inversion (\ref{equation: at the identity}), we see that:

$$\sum_{\nu \in \eigs} \hat{k}(\nu) h(t \nu_{\infty}) \geq
\vol(\tilde{\Omega}') \int \hat{k}(\nu) h(t \nu_{\infty}) d \planch(\nu)
- c |Z| t ^{-d+1/2}, \ \ \ (0 < t \leq 1).  $$

It follows:
$$\liminf_{t \rightarrow 0} t ^{d} \sum_{\nu \in \eigs}
h(t \nu_{\infty})
\geq \vol(\tilde{\Omega}') \liminf_{t \rightarrow 0} t ^ d \frac{\int \hat{k}(\nu) h(t \nu_{\infty})
d \planch(\nu) }{\| k \|_{\specnorm}}.$$

We apply Lem. \ref{lem:main} to choose $k$, taking $k = k_n$ for some sufficiently large $n$
and then take $n \rightarrow \infty$. 
Note moreover
that we can choose $\vol(\tilde{\Omega}')$ to be arbitrarily close to  $\vol(\tilde{\Omega})=\vol(\Omega)$,
which in turn (by choosing $\Omega$ sufficiently large) can be made arbitrarily close
to $\vol(\Gamma \backslash G_S)$.  We see that:
\begin{equation}
\label{eq:almostdone}
\liminf_{t \rightarrow 0} t ^{d} \sum_{\nu \in \eigs}
h(t \nu_{\infty}) \geq \vol(\Gamma \backslash G_S)
\liminf_{t \rightarrow 0 }\left( t ^{d} \int h(t \nu_{\infty}) d \planch(\nu)\right).
\end{equation}

We choose $H^{\infty}, h$ as in Lem. \ref{lem:hchoice}; this choice depends
on a parameter $\epsilon$, which we then  let approach $0$. 
The same reasoning as used in Sec. \ref{sec:Weyllawcompact} (see esp.
  \eqref{tail estimate} and \eqref{eq:spongebob})
 allows one to deduce from (\ref{eq:almostdone}) a lower bound in the Weyl law; in the
 present context,  we  
must use as an {\em a priori} input Donnelly's upper bound
(\ref{eq:donnelly}). In any case, we obtain:
\begin{equation} \label{eq:lowerbound} \liminf_{T \rightarrow \infty}
\frac{\#\{ \nu \in \eigs: \langle \nu_{\infty}, \nu_{\infty} \rangle \leq T \}}{T ^{d/2}}
\geq \vol(\Gamma \backslash G_S) \alpha(\G).\end{equation} 
This is the lower bound in the Weyl law. (It can be verified by computation
that the lower bound is ``correct'', i.e., the constant coincides with that
described at the start of the Introduction; we sketch a different proof of this correctness in \S \ref{subsec:final}.)
\qed


%

%

%

\subsection{Proof of the main lemma} \label{subsec:construc}

We now turn to the proof of Lem. \ref{lem:hchoice}. 
We shall first prove a preparatory result which asserts, in essence, that linear combinations of
trigonometric functions do not take small values too often.

\begin{lemma} \label{lem:sillier}
Let $A= \mathbb{R}^ d \times
(\mathbb{R}/\mathbb{Z})^{e}$, and let $B(T)$ be the
subset of $A$ consisting of elements whose projection to $\mathbb{R}^ d$
lies in a Euclidean ball of radius $T$ and centered at $0$.
Let $\chi_1, \dots, \chi_n$ be distinct unitary characters of
$A$, let $0 \neq a_i \in \mathbb{C}$, and put
$F = \sum_{i=1}^{n} a_i \chi_i$.  Fix a Haar measure $d \mu$ on $A$.

Then:
\begin{equation}\label{eq:lax}
\lim_{\epsilon \rightarrow 0}
\left(
\lim_{T \rightarrow \infty} \frac{
\mu \{ x \in B(T): |F(x)| \leq \epsilon \}}{T ^ d} \right) = 0.\end{equation}
\end{lemma}

\proof
The map $\underline{\chi} = (\chi_1, \dots, \chi_n)$
maps $A$ to the torus $\mathbb{T}  :=
\{ (z_1, \dots, z_n) \in \mathbb{C}^ n: |z_i| = 1 \text{ for all } i \}$.
Let $L$ be the function $(z_1, \dots, z_n) \rightarrow \sum_{i} a_i z_i$
on $\mathbb{T}$, and let $\mathbb{T}'$ be the closure
of the image $\underline{\chi}(A)$. Then $\mathbb{T}'$ is a subtorus of $\mathbb{T}$. 
$L$ does not vanish identically on $\mathbb{T}'$,
since otherwise $F = 0$, contradicting the linear independence of
characters.

Let $d \mu_T$ be the measure $\frac{d \mu}{\mu(B(T))}$ restricted to $B(T)$,
and let $\nu_T = \underline{\chi}_{*} \mu_T$.
Any weak limit (as $T \rightarrow \infty$) of the measures $\nu_T$ is supported
on $\mathbb{T}'$ and invariant by $\mathbb{T}'$; consequently,
the $\nu_T$ converge to the invariant probability measure $\nu$
on $\mathbb{T}'$.
(\ref{eq:lax}) amounts to the assertion that the
zero-locus of $L$ on $\mathbb{T}'$ has zero measure w.r.t $\nu$; but this is
obvious, as $L$ is real-analytic and nonvanishing.
\qed

\proof (of Lem. \ref{lem:main})
For fixed $g$ s.t. $g_{\infty} \notin K_{\infty} $,
the spherical functions
$\Xi_{\nu}(g)$ decay as $\nu_{\infty} \rightarrow \infty$. An elegant quantification of  this
has been given by Duistermaat, Kolk and Varadarajan in \cite[Cor. 11.2]{DKV};
it shows in particular that for $g_{\infty}$ in a fixed compact
subset excluding $K_{\infty}$,
we have a bound of the form $|\Xi_{\nu}(g)| \leq c' (1+\| \nu_{\infty}\|)^{-1/2}$
for some $c'$. Combining with Plancherel inversion
(\ref{eq:planchinversion}) and  (\ref{eq:planchandhaar}), it follows that there is $c''$ such that
\begin{multline}
\sup_{g \in F} t ^{d-1/2} |k_t(g)| \\
\leq c'' t ^{d-1/2} \| k \|_{\specnorm} \int_{\nu_{\infty} \in \mathfrak{a}_{\infty}^{*}}
h(t \nu_{\infty}) (1+\| \nu_{\infty}\|)^{d-r-1/2} d \haar(\nu_{\infty})
\\ \leq c'' \| k \|_{\specnorm} \int_{\mathfrak{a}_{\infty}^*}
h(\nu_{\infty}) (t+\| \nu_{\infty}\|)^{d-r-1/2}
d \haar(\nu_{\infty}),\end{multline}
whence (\ref{eq:kestclose}).

Next, we turn to the construction of the $k_n$s. Apply Lem. \ref{lem2} with $A = \mathfrak{a}_S$, $W= W_S$,
$T = \{ \Gamma_{A,i}: 1 \leq i \leq I.\}$
 Lem. \ref{lem2} produces
a nonzero, compactly supported, $W_S$-invariant distribution $f$ on $\mathfrak{a}_S$
which vanishes on any $\Gamma_{A,i}$-invariant function, for $1 \leq i \leq
I$; moreover, this $f$ is a finite linear combination of point masses. 
Let $\check{f}(x) = \overline{f(-x)}$, and
replace $f$ by $f \star \check{f}$; then $f$
is still $W_S$-invariant and a linear combination of point masses,
and $\hat{f}(\nu)$ is a non-negative real number whenever $\nu,
\overline{\nu}$ are $W_S$-conjugate.

The function $\nu \mapsto \hat{f}(\nu)$ is a finite linear
combination of characters.
Put $K = \sup_{\nu \in \mathfrak{a}_{S}^{*}: \pi(\nu) \,
\mathrm{unitary}}\hat{f}(\nu).$
The supremum is easily verified to be finite.
Now set $P_n(x) = 1- (1-\frac{x ^ 2}{K ^ 2})^ n$.
Then $P_n$ satisfies:
\begin{enumerate}
\item For each $n$,
$P_n(0) = 0$.
\item For each $n$, $0 \leq P_n(x) \leq 1$ for $x \in [-K,K]$.
\item The $P_n$ converge uniformly on any compact subset
of $[-K,K] - \{ 0 \}$ to the constant function with value $1$.
\end{enumerate}

Equip the space of compactly supported distributions
on $\mathfrak{a}_{S}$ with the algebra structure in which
multiplication corresponds to convolution.
With this in mind, set
$f_{n} = P_{n}(f)$; in other words,
$\widehat{f_{n}}(\nu) = P_{n}(\hat{f}(\nu))$. 
Thus $f_n$ is a $W_S$-invariant distribution on $\mathfrak{a}_{S}$, which
is indeed a finite linear combination of point masses. 

We now choose $k_n \in \mathcal{D}(K_S \backslash G_S / K_S)$ such that $\Satake k_n = f_n$; or equivalently
(by \eqref{eqn:cd}) $\widehat{k_n}(\nu) = \widehat{f_n}(\nu)$.
(Here we have used a distributional extension of the Satake isomorphism;
however, as we have remarked, one can avoid this entirely
and use Thm. \ref{theorem: Paley-Wiener} as stated, 
because $k_n$ enters only through its convolution with an appropriate $H_t$).
It is clear that $k_n$ satisfy (\ref{mainfirst}) and (\ref{mainsecond}) of Lem. \ref{lem:main}. 
Moreover, by choice, $\|k_n\|_{\specnorm} \leq 1$. 

It remains to verify (\ref{mainthird}) of Lem. \ref{lem:main}. For this it suffices to check
that, for $h$ a non-negative Schwarz function on $\mathfrak{a}_{\infty, \temp}^{*}$, we have:
\begin{equation} \label{eq: de}\lim_{n \rightarrow \infty} \left(\liminf_{t \rightarrow 0}
\frac{\int h(t \nu_{\infty})\widehat{f_{n}}(\nu)
d \planch(\nu)}{\int h(t \nu_{\infty}) d \planch(\nu)} \right) =   1.\end{equation}

It is easy to see that $\int
h(t \nu_{\infty}) d \planch(\nu) \sim c t ^{-d}$, for some $c > 0$.
To check (\ref{eq: de}), it suffices to check that the limit
$$L := \lim_{\epsilon \rightarrow 0} \left( \limsup_{t \rightarrow 0}
t ^{d} \int_{\nu
 \in
\mathfrak{a}_{S,\temp}^{*}: |\hat{f}(\nu)| \leq \epsilon} h(t \nu_{\infty})
d \planch(\nu)\right)$$ equals $0$.
From (\ref{eq:planchandhaar}),
we deduce that
\begin{multline} \label{eq:almostpre}
L \leq  c_1
\lim_{\epsilon \rightarrow 0} \left(
\limsup_{t \rightarrow 0} t^d\int_{\nu \in \mathfrak{a}_{S,\temp}^{*}: |\hat{f}(\nu)| \leq \epsilon}
(1 + \|\nu_{\infty} \|)^{d-r}| h(t \nu_{\infty})| d \haar(\nu)\right).
\end{multline}

Let $c_1'$ be so that $\sup_{\nu \in \mathfrak{a}_{\infty,\temp}^{*}}
|h(\nu)| (1+\| \nu \|)^{d+1} \leq c_1'$. Then:
\begin{multline}\label{eq:almostpr}
L  \leq
c_1 c_1' \lim_{\epsilon \rightarrow 0} \limsup_{t \rightarrow 0}
t^r \int_{\nu \in
\mathfrak{a}_{S,\temp}^{*}: |\hat{f}(\nu)| \leq \epsilon}
\frac{(t+\|t \nu_{\infty} \|)^{d-r} }{(1+\|t \nu_{\infty}\|)^{d+1}} d \haar(\nu)\\
\leq 
c_1 c_1' \lim_{\epsilon \rightarrow 0} \limsup_{t \rightarrow 0} t^r
\int_{\nu \in
\mathfrak{a}_{S,\temp}^{*}: |\hat{f}(\nu)| \leq \epsilon}
(1+\|t \nu_{\infty}\|)^{-r-1} d \haar(\nu)
\end{multline}
From Lem. \ref{lem:sillier}, we conclude that the quantity
on the right of (\ref{eq:almostpr}) is $0$, i.e. $L=0$ as claimed.
This completes the verification of condition (\ref{mainthird}) of the Lemma. 
\qed
\subsection{Deduction of Thm. \ref{thm:mainthm} from the $S$-arithmetic Weyl law.} \label{subsec:final} 
We now discuss the final reductions to complete the proof of Thm. \ref{thm:mainthm}. The argument, in words, is as follows:  (\ref{eq:lowerbound}) asserts, in effect, the lower bound in Weyl law for the cuspidal
spectrum on a certain finite union of locally symmetric spaces. Since one knows, by Donnelly's work,
the upper bound in the Weyl law for each of these spaces, the Weyl law for each individual
locally symmetric space follows.


The number of $(\Gamma, G_{\infty} K_S)$ double cosets in $G_S$ is finite;
in particular, there is a finite collection $(g_i)$ such that $G_S$
is the disjoint union of the double cosets $\Gamma g_i G_{\infty} K_S$.
Without loss of generality we may take $g_1 = 1$.

Put $\Gamma_{\infty, i} =g_i^{-1} \Gamma  g_i \cap G_{\infty} K_S $, which
we regard (by projection) as a lattice in $G_{\infty}$. Then the space
$\Gamma \backslash G_S /K_S$ is the disjoint union
of the locally symmetric spaces $\Gamma_{\infty,i} \backslash G_{\infty} /K_{\infty}$.
Let $N(T)$ (resp. $N_i(T))$ be the number of cuspidal eigenfunctions of the Laplacian
on $\Gamma \backslash G_S/K_S$ (resp. $\Gamma_{\infty,i} \backslash G_{\infty} /K_{\infty}$) 
with eigenvalue $\leq T$. 
Then $\sum_{i} N_i(T) = N(T)$.

Now (\ref{eq:lowerbound}) shows that $\liminf_{T \rightarrow \infty} \frac{N(T)}{T ^{d/2}}
\geq \alpha(\G) \vol(\Gamma \backslash G_S)$.
On the other hand, the main result of \cite{donnelly} asserts that
\begin{equation} \label{eq:upperbound} \limsup_{T \rightarrow \infty} \frac{N_i(T)}{T ^{d/2}} \leq c(\Gamma_{\infty,i} \backslash G_{\infty} /K_{\infty}),\end{equation} where
$c(M)$ is defined in the first paragraph of the Introduction.

We now claim that
\begin{equation} \label{eq:equality} \sum_{i} c(\Gamma_{\infty,i} \backslash G_{\infty} /K_{\infty})
= \alpha(\G) \vol(\Gamma \backslash G_S).\end{equation} Actually, this can be done using the explicit
forms of the Plancherel measure, etc.; however, let us sketch a slightly more ``conceptual''
approach, which has the disadvantage of invoking slightly more.
It is a theorem of A. Borel \cite{armand} that $G_{\infty} = \G(\R)$ admits a cocompact torsion-free arithmetic lattice $\Lambda$.
By the discussion of Sec. \ref{sec:Weyllawcompact} -- applied in the case $S = \{ \infty \}$ -- the number of Laplacian
eigenfunctions on $\Lambda \backslash G_{\infty} / K_{\infty}$ with
eigenvalue $\leq T$ is $\sim \alpha(\G) \vol(\Lambda \backslash G_{\infty}) T ^{d/2} $.  On the other hand, Weyl's original proof applies to the compact Riemannian
manifold $\Lambda \backslash G_{\infty} / K_{\infty}$ and shows that
this number is also $\sim c(\Lambda \backslash G_{\infty} / K_{\infty}) T ^{d/2}$; consequently
we have the equality $\alpha(\G) \vol(\Lambda \backslash G_{\infty}) =
c(\Lambda \backslash G_{\infty}/K_{\infty})$. Now it is easy to see that
both sides of this equality change by the same multiplicative factor when
$\Lambda$ is replaced by any other (not necessarily cocompact) lattice $\Lambda'$; in particular,
$\alpha(\G) \vol(\Gamma_{\infty,i} \backslash G_{\infty}) = c(\Gamma_{\infty,i} \backslash G_{\infty}/K_{\infty})$; summing over $i$, and noting that $\vol(\Gamma \backslash G_S) =
\sum_{i} \vol(\Gamma_{\infty,i} \backslash G_{\infty})$, we have proven (\ref{eq:equality}).

From (\ref{eq:equality}), (\ref{eq:lowerbound}) and (\ref{eq:upperbound}), it follows that
for each $i$ we have $\lim_{T \rightarrow \infty} \frac{N_i(T)}{ T ^{d/2}} = c(\Gamma_{\infty,i} \backslash G_{\infty} /K_{\infty})$, i.e. Weyl's law holds for the symmetric space $\Gamma_{\infty,i} \backslash G_{\infty}/K_{\infty}$.
On the other hand, any congruence subgroup of $\G(\Z)$ arises as $\Gamma_{\infty,i}$
for some $S$-arithmetic congruence subgroup $\Gamma \subset \G(\Z[S ^{-1}])$ and
some $i$. Theorem \ref{thm:mainthm} follows. \qed

\appendix
\section{ Existence of cusp forms for general $\G$ --- a proof using Whittaker functions}

The considerations of the previous sections give the full Weyl law. However,
we would also like to explain a short proof of the {\em existence} of cusp 
forms that uses Whittaker functions; this does not give the Weyl law but,
unlike the proof of the Weyl law, gives a very explicit method for constructing
cuspidal functions. 
It is easy to modify the present proof to show (e.g.)
the existence of infinitely many cusp forms which do not
lie in the image of a fixed set of functorial lifts.


For $T \geq 0$, set $A_{S,> T} = \{ a \in A_S:
|\alpha(a)|_{S} > T,  \ \forall \ \alpha \in \Delta \}$.
Define $\Siegel_{T} \subset G_S$ via $\Siegel_{ T} = N_S \cdot A_{S, > T} \cdot K_S$
and set $\Gamma_{B} = B_S \cap \Gamma$, $\Gamma_{A} = A_S \cap \Gamma$,
$\Gamma_N = N_S \cap \Gamma$.  Then $\Gamma_B . \Siegel_T \subset \Siegel_T$.
The image of $\Siegel_T$ in $\Gamma \backslash G_S$ is called a ``Siegel set.''

We shall need the following
properties:

\begin{enumerate}
\item For every $T$, (the image of) $\Siegel_{T}$ intersects every connected component
of $G_S/K_S$.

(Let $\mathfrak{t} = \mathrm{Lie}(\A)$, so we have an exponential
map $\exp: \mathfrak{t}(\R) \rightarrow \A(\R) \subset \G(\R)$.
Let $a \in \mathfrak{t}$ belong to the interior of the positive Weyl
chamber. Then for any $x \in G_S/K_S$ and sufficiently large $t > 0$,
$\exp(t . a) x$ belongs to $\Siegel_T$.)

\item \label{reductiontheory} If $T$ is chosen sufficiently large,
the natural map $$\iota_{T}: \Gamma_B \backslash \Siegel_{T}
\rightarrow \Gamma \backslash G_S$$ is a homeomorphism
onto an open subset $\mathscr{S}_T$ of $\Gamma \backslash G_S$.

(This is $S$-arithmetic reduction theory.)
\item If $Z \subset G_S $ is any compact subset there exists
$T' > T$ so that $\Siegel_{T'} . Z \subset \Siegel_T$.

(Set $Z' = K_S. Z$. Since $Z'$ is compact, there exists $L \geq 1$
such that $Z' \subset N_S \cdot \{ a \in A_S, L ^{-1} < |\alpha(a)|_S < L
, \forall \alpha \in \Delta \} \cdot K_S$. Then $T' = LT$ suffices.)

\end{enumerate}

\begin{prop} \label{prop2}
Suppose $k \in C ^{\infty}_c(K_S \backslash G_S/K_S)$ is nonzero.
Then
the convolution operator $f \mapsto f \star k$ is nontrivial
on $L ^ 1_{\mathrm{loc}}(\Gamma \backslash G_S)$.

\end{prop}

\proof

We construct a function $f$ on $\Gamma_{B}\backslash
G_S$ such that $f \star k = \lambda f$ for some $\lambda \neq 0$,
and so that $f|\Siegel_T \neq 0$ for all $T$.

This immediately implies the assertion of Proposition
\ref{prop2} by ``transporting'' $f$
to $\Gamma \backslash G_S$. More formally: let $T$
be sufficiently large. Then $f \circ \iota ^{-1}$
defines a function on $\mathscr{S}_T$. Extending this
function by $0$ off $\mathscr{S}_T$, we obtain a function $f'
\in L ^ 1_{\mathrm{loc}}(\Gamma_S \backslash G_S)$.
We claim $f' \star k \neq 0$.
Indeed, let $T' > T$ be so that $\Siegel_{T'} . \mathrm{supp}(k)
\subset \Siegel_{T}$. Then for $x \in \Siegel_{T'}$
we have $f' \star k(\iota(x)) = f \star k(x) = \lambda f(x)$.
Since $\lambda \neq 0$ and $f| \Siegel_{T'} \neq 0$, we are done.

Let $\psi$ be a nondegenerate character of $N_S$ trivial on
$\Gamma_N$.
(Nondegenerate means its stabilizer in $A_S$ is central in $G_S$,
i.e. trivial by the adjointness assumption.) Let $\pi_S$ be a generic spherical representation of $G_S$
upon which $k$ acts nontrivially. (It exists, since the spherical
constituent of $\pi(\nu)$ is generic for an open
dense set of $\nu \in \mathfrak{a}_S ^{*}$
(\cite{jacquet}, \cite{rodier})
and the transform
$\hat{k}(\nu)$, being holomorphic, cannot vanish identically
on such a set.)
Let $W$ be a spherical Whittaker function corresponding to
$\pi_S$, i.e. $W$ is a function on $G_S$ of moderate
growth whose right translates realize the representation
$\pi_S$, and satisfying $W(ng) = \psi(n) W(g)$.
Then $W \star k = \lambda W$ for some $\lambda \neq 0$.

Set
\begin{equation} \label{eqn:f}
f(g) = \sum_{\tau \in \Gamma_N \backslash \Gamma_B} W(\tau g).\end{equation}
The series converges absolutely and uniformly
in compacta; this may be deduced from the ``rapid decay'' of $W$.

%

$f$ defines a function on $G_S$ that is left-invariant
by $\Gamma_B$ and right-invariant by $K_S$.
To check that $f | \Siegel_{T}$ is nonvanishing
it suffices to check that the function $g \mapsto \int_{n \in
\Gamma_N \backslash N_S} f(ng) \psi(n) dn$ is not identically
zero for $g \in \Siegel_T$.
Note that each term $W(\tau g)$ in the definition (\ref{eqn:f})
of $f$ transforms (on the left)
by the character $n \mapsto \psi(\tau n \tau ^{-1})$ of $N_S$,
and these characters are all distinct since $\psi$ is nondegenerate.
Thus $\int_{n \in \Gamma_N \backslash N_S} f(ng) \psi(n) dn = W(g)$
and it suffices to check that $W(g)|\Siegel_T \neq 0$ for any $T$.

However, note that $G_S/K_S$ is a countable
union of real symmetric spaces, and $W(g)$ defines an eigenfunction
of the Laplacian on each component; in particular, it is real-analytic
on each component. It follows that there is at least one component
such that $W(g)$ does not vanish identically on any open set.

Since the image of $\Siegel_T$
intersects every connected component of $G_S/K_S$ it
follows that $f|\Siegel_{T}$ is nonzero, and we are done. \qed

That there exists at least one cusp form now follows from 
 Prop. \ref{prop1} and Prop. \ref{prop2} 
 and Cor. \ref{cor1}. It is simple to modify this argument to show, e.g., that
 there are infinitely many cusp forms, or infinitely many cusp forms which
 are not self-dual. 


\begin{thebibliography}{1}

\bibitem{armand}
A. Borel,
\newblock Compact Clifford-Klein forms of symmetric spaces.
\newblock {\em Topology} 2:111-122, 1963.

\bibitem{cartier}
P. Cartier,
\newblock Representations of $p$-adic groups: a survey.\newblock
{\em Automorphic forms, representations and {$L$}-functions}.
\newblock American Mathematical Society, Providence, RI, 1979.

\bibitem{donnelly}
Harold Donnelly,
\newblock On the cuspidal spectrum for finite volume symmetric spaces.
\newblock {\em J. Differential Geom.} 17:239--253, 1982.

\bibitem{DKV}
J. Duistermaat, J. Kolk and V. Varadarajan.
\newblock
Functions, flows and oscillatory integrals on flag manifolds and
conjugacy classes in real semisimple Lie groups.
\newblock      {\em    Compositio Math.} 49:309--398, 1983.

\bibitem{DKV2}
J. Duistermaat, J. Kolk and V. Varadarajan.
\newblock
Spectra of compact locally symmetric manifolds of negative curvature. 
\newblock      {\em   Invent. Math.} 52:27--93, 1979.

\bibitem{FK}
Y. Flicker and D. Kazhdan.
\newblock
A simple trace formula.
\newblock {\em J. Analyse Math.}, 50: 189--200, 1988. 


\bibitem{gangolli}

R.~Gangolli.
\newblock On the {P}lancherel formula and the {P}aley-{W}iener theorem for
spherical functions on semisimple {L}ie groups.
\newblock {\em Ann. of Math. (2)}, 93:150--165, 1971.

\bibitem{helgason}
Sigurdur Helgason.
\newblock An analogue of the {P}aley-{W}iener theorem for the {F}ourier
transform on certain symmetric spaces.
\newblock {\em Math. Ann.}, 165:297--308, 1966.


\bibitem{Helgason-topics}
Sigurdur Helgason.
\newblock {\em Topics in harmonic analysis on homogeneous spaces}, volume~13 of
{\em Progress in Mathematics}.
\newblock Birkh\"auser Boston, Mass., 1981.

\bibitem{gga}
S. Helgason. \newblock
Groups and geometric analysis.
\newblock American Mathematical Society, Providence, RI, 2000.
\newblock Mathematical Surveys and Monographs, No. 83.

\bibitem{iwaniec}
Henryk Iwaniec.
\newblock {\em Spectral methods of automorphic forms}, volume~53 of {\em
Graduate Studies in Mathematics}.
\newblock American Mathematical Society, Providence, RI, second edition, 2002.






\bibitem{jacquet}
Herv{\'e} Jacquet.
\newblock Fonctions de {W}hittaker associ\'ees aux groupes de {C}hevalley.
\newblock {\em Bull. Soc. Math. France}, 95:243--309, 1967.

\bibitem{labessemuller}
J.-P. Labesse and W. M{\"u}ller.
\newblock Weak Weyl's law for congruence subgroups.
\newblock {\tt math.RT/0404037}.



\bibitem{Lax-Phillips-automorphic-book}
Peter~D.  Lax and Ralph~S.  Phillips.
\newblock {\em Scattering theory for automorphic functions}.  \newblock
Princeton Univ.  Press, Princeton, N.J., 1976.
\newblock Annals of Mathematics Studies, No.  87.

\bibitem{macdonald}
I. G. Macdonald.
\newblock Spherical functions on a group of $p$-adic type.
\newblock Publications of the Ramanujan Institute, volume 2, 1971.

\bibitem{miller}
Stephen~D. Miller.
\newblock On the existence and temperedness of cusp forms for {${\rm SL}\sb
3({\mathbb Z})$}.
\newblock {\em J. Reine Angew. Math.}, 533:127--169, 2001.

\bibitem{mw}
C. Moeglin and J.-L. Waldspurger,
\newblock D\'ecomposition spectrale et s\'eries d'Eisenstein.
\newblock {\em Progress in Mathematics}, 113, Birkh{\"a}user, 1994.

\bibitem{muller}
Werner M{\"u}ller.
\newblock Weyl's law for the cuspidal spectrum of {${\rm SL}\sb n$}.
\newblock {\em C. R. Math. Acad. Sci. Paris}, 338(5):347--352, 2004.

\bibitem{ps}
R.~Phillips and P.~Sarnak.
\newblock Automorphic spectrum and {F}ermi's golden rule.
\newblock {\em J. Anal. Math.}, 59:179--187, 1992.
\newblock Festschrift on the occasion of the 70th birthday of Shmuel Agmon.



\bibitem{rodier}
Fran{\c{c}}ois Rodier.
\newblock Whittaker models for admissible representations of reductive
{$p$}-adic split groups.
\newblock In {\em Harmonic analysis on homogeneous spaces (Proc. Sympos. Pure
Math., Vol. XXVI, Williams Coll., Williamstown, Mass., 1972)}, pages
425--430. Amer. Math. Soc., Providence, R.I., 1973.

\bibitem{Sarnak-cusp-II}
P.~Sarnak.
\newblock On cusp forms. {II}.
\newblock In {\em Festschrift in honor of I. I. Piatetski-Shapiro on the
  occasion of his sixtieth birthday, Part II (Ramat Aviv, 1989)}, volume~3 of
  {\em Israel Math. Conf. Proc.}, pages 237--250. Weizmann, Jerusalem, 1990.

\bibitem{selberg}
A.~Selberg.
\newblock Harmonic analysis and discontinuous groups in weakly symmetric
{R}iemannian spaces with applications to {D}irichlet series.
\newblock {\em J. Indian Math. Soc. (N.S.)}, 20:47--87, 1956.


\bibitem{Tits-Corvallis}
J. Tits. 
\newblock Reductive groups over local fields. 
\newblock
{\em Automorphic forms, representations and {$L$}-functions}.
\newblock American Mathematical Society, Providence, RI, 1979.


\bibitem{warner}
Garth Warner.
\newblock {\em Harmonic analysis on semi-simple {L}ie groups. {II}}.
\newblock Springer-Verlag, New York, 1972.
\newblock Die Grundlehren der mathematischen Wissenschaften, Band 189.

\bibitem{wolpert}
Scott Wolpert. 
\newblock{Disappearance of cusp forms in special families.}
\newblock {\em Ann. of Math. (2)}, 139:239--291, 1994. 

\end{thebibliography}
\end{document}